%% file: HISP.tex
\documentclass[journal]{IEEEtran}

\usepackage[xindy]{glossaries}
\usepackage[hidelinks]{hyperref}
\usepackage{amsfonts,amssymb,mathtools,bm,amsthm,xfrac}
\usepackage[OT1]{fontenc}
\usepackage[inline]{enumitem}
\usepackage{graphicx,caption,subcaption,xcolor}
\usepackage{cleveref, siunitx, algorithmic, algorithm}

\newtheorem{theorem}{Theorem}

\newtheorem{lemma}{Lemma}

\theoremstyle{definition}

\newtheorem{example}{Example}
\newtheorem{remark}{Remark}

\def\eqns#1{\begin{equation*}#1\end{equation*}}
\def\eqnsml#1{\begin{multline*}#1\end{multline*}}
\def\eqnl#1#2{\begin{equation}\label{#1}#2\end{equation}}
\def\eqnla#1#2{\begin{align}\label{#1}#2\end{align}}
\def\eqnsa#1{\begin{align*}#1\end{align*}}

\def\zero{\mathbf{0}}
\def\one{\mathbf{1}}
\def\bsF{\bm{F}}
\def\bsi{\bm{i}}
\def\bsj{\bm{j}}
\def\bsk{\bm{k}}
\def\bfm{\mathbf{m}}
\def\bso{\bm{o}}
\def\bsP{\bm{P}}
\def\bsV{\bm{V}}
\def\boX{\mathbf{X}}
\def\bfx{\mathbf{x}}
\def\boZ{\mathbf{Z}}
\def\bfz{\mathbf{z}}
\def\bsphi{\bm{\phi}}	

\def\calA{\mathcal{A}}
\def\calN{\mathcal{N}}
\def\calO{\mathcal{O}}
\def\calP{\mathcal{P}}

\def\bbI{\mathbb{I}}
\def\bbJ{\mathbb{J}}
\def\bbL{\mathbb{L}}
\def\bbN{\mathbb{N}}
\def\bbO{\mathbb{O}}

\def\bbR{\mathbb{R}}
\def\bbT{\mathbb{T}}

\def\c{\mathrm{c}}
\def\d{\mathrm{d}}
\def\m{\mathrm{m}}
\def\u{\mathrm{u}}

\def\ex{\mathrm{ex}}

\DeclareMathOperator{\graph}{Gr}
\DeclareMathOperator*{\argmax}{argmax}

\def\spacedAnd{\quad\mbox{ and }\quad}
\def\defeq{\doteq}
\def\given{\,|\,}
\def\ind#1{\one_{#1}}
\def\sst{\,:\,} 
\def\st{:}

\def\inx{\bullet}
\def\op{\flat}
\def\obj{\sharp}

\def\Pr{\mathrm{Pr}}
\def\bp{\bar{p}}
\def\bq{\bar{q}}

\newacronym{hisp}{HISP filter}{hypothesised filter for independent stochastic population}
\newacronym{disp}{DISP filter}{filter for distinguishable and independent stochastic populations}
\newacronym{phd}{PHD}{probability hypothesis density}
\newacronym{lmb}{LMB}{Labelled multi-Bernoulli}
\newacronym{cphd}{CPHD}{cardinalised PHD}
\newacronym{mht}{MHT}{multiple hypothesis tracking}
\newacronym{mtt}{MTT}{multi-target tracking}
\newacronym{mc}{MC}{Monte Carlo}
\newacronym{smc}{SMC}{Sequential Monte Carlo}
\newacronym{kf}{KF}{Kalman filter}
\newacronym{rfs}{RFS}{random finite set}
\newacronym{snr}{SNR}{signal-to-noise ratio}

\title{Multi-target filtering with linearised complexity}

\author{Jeremie Houssineau and Daniel E. Clark%
\thanks{J.\ Houssineau is with the Department of Statistics and Applied Probability, National University of Singapore, 117546, SG. Email: stahje@nus.edu.sg}%
\thanks{D.E. Clark is with the School of Electrical and Physical Sciences, Heriot-Watt University, Edinburgh EH14 4AS, UK. Email: D.E.Clark@hw.ac.uk.}}

\begin{document}


\maketitle

\begin{abstract}
An algorithm for the estimation of multiple targets from partial and corrupted observations is introduced based on the concept of partially-distinguishable multi-target system. It combines the advantages of engineering solutions like \gls{mht} with the rigour of point-process-based methods. It is demonstrated that under intuitive assumptions and approximations, the complexity of the proposed multi-target estimation algorithm can be made linear in the number of tracks and linear in the number of observations, while naturally preserving distinct tracks for detected targets, unlike point-process-based methods.
\end{abstract}



\section{Introduction}

\IEEEPARstart{M}{ulti}-target tracking refers to the problem of estimating the number and state of targets in a dynamical multi-target system via partial and uncertain observations corrupted by detection failures and false alarms. The difficulty of this task is further amplified by the fact that targets in the system appear and disappear at unknown times. The main challenge with \gls{mtt} is the absence of a priori information about \emph{data association}, i.e.\ about the association between the received observations and the different targets being estimated. This aspect makes the problem highly combinatorial in nature so that approximations are required for systems made of a large number of targets and/or for which the data association problem is not easily resolved.

The existing solutions in the field of \gls{mtt} can be divided into two classes of methods. One class consists of ``classical'' methods \cite{Blackman1986,BarShalom1987,Okuma2004,Vermaak2005} that are based on practical generalisations of single-target estimation. Their strength lies in their ability to distinguish the targets and to naturally characterise each of them. The other class of methods comprises approaches based on simple point processes, such as the \gls{phd} filter introduced by \cite{Mahler2003} or one of the many variants of it \cite{Singh2009,Vo2009,Pace2013}. These methods successfully propagate global statistics about the system of interest and integrate false alarms and appearance of targets in a principled way. However, they do not naturally propagate specific information about targets because of the point process assumption of indistinguishability. One of the attempts to overcome this limitation \cite{Vo2013,Vo2014} uses marked point processes.

More recently, a new paradigm for modelling systems of targets with uncertain cardinality and state has been introduced in \cite{Houssineau2016}, which has been combined in \cite{Houssineau2015} with a novel representation of uncertainty \cite{Houssineau2017,Houssineau2016_dataAssimilation} allowing for modelling partial information. This paradigm embeds the concept of partially-distinguishable systems which is useful in \gls{mtt}, e.g.\ for jointly representing never-detected targets (indistinguishable) and previously-detected ones (distinguishable). It also enables the computation of more diverse types of statistics \cite{Houssineau2016} about multi-target systems than with point processes \cite{Delande2014_Var}, which are helpful, for instance, for sensor management \cite{Delande2014_SensorControl}. By considering the usual assumptions of \gls{mtt}, a first algorithm following from this paradigm has been introduced in \cite[Chapt.\ 3]{Houssineau2015} and detailed in \cite{Delande2016,Delande2017_SSA} and is referred to as the \gls{disp}. This filter can be derived without relying on approximations but unsurprisingly displays high computational complexity, similarly to the \gls{mht}. In this article, an additional \gls{mtt} algorithm, called the \gls{hisp}, is derived and its efficiency is demonstrated on simulated data. The objective is to have recourse to some intuitive approximations in order to lower the complexity.

The structure of the article is as follows: the recursion of the \gls{hisp} is given in \cref{sec:HISP:derivation} and in a more practical way in \cref{sec:HISP:pdf}. This is followed by an introduction of the main approximation in \cref{sec:mainApproximation}, showing how the computational complexity is brought down to linear in the number of tracks and in the number of observations. A discussion of the connections with existing algorithms can be found in \cref{sec:relationWithOtherWorks}. Details of the implementation are given in \cref{sec:implementation} and performance is demonstrated on simulated data in \cref{sec:simulations}.

\section[The HISP filter]{The \gls{hisp}}
\label{sec:HISP:derivation}

The specific modelling enabling the derivation of the \gls{hisp} is given in \cref{ssec:individualModelling} and \cref{ssec:populationModelling}, followed by a detailed presentation of the successive steps of the \gls{hisp}'s recursion in \cref{ssec:timeFiltering} and \cref{ssec:observationFiltering}.

Without loss of generality, the time is indexed by the set $\bbT \defeq \bbN$, where $\bbN$ is the set of natural numbers. For any $t \in \bbT$, the state and observation spaces, denoted $\boX = \boX^{\inx} \cup \{\psi_{\obj},\psi_{\op}\}$ and $\boZ = \boZ^{\inx} \cup \{\phi\}$ respectively, are defined as the union of the spaces of interest $\boX^{\inx}$ and $\boZ^{\inx}$ together with the empty states $\psi_{\obj}$ and $\psi_{\op}$ as well as the empty observation $\phi$. It is common to consider $\boX^{\inx} = \bbR^d$ and $\boZ^{\inx} = \bbR^{d'}$ with the integers $d$ and $d'$ usually verifying $d > d'$. The empty state $\psi_{\obj}$ is attributed to targets that are not in the area of interest, modelled by $\boX^{\inx}$ at time $t$, but which might appear in it at any subsequent time. The state $\psi_{\op}$ is attributed to \emph{false-alarm generators} which are objects in the field of view of the sensor that are not of direct interest but which might interfere with the observation of the targets\footnote{this construction has been previously proposed in \cite{Singh2009}}. Note that the symbol ``$\obj$'' will generally refer to targets whereas ``$\op$'' will refer to the false alarm. The empty observation $\phi$ is attributed to objects that are not actually observed at a given observation time. The use of empty state/observation is common in the target-tracking literature \cite{Caron2011_conditional}. Because of the nature of the spaces $\boX$ and $\boZ$, the integral $\int f(\bfx) \d\bfx$ of an integrable function $f$ over $\boX$ is understood as
\eqns{
\int f(\bfx) \d\bfx = \int_{\boX^{\inx}} f(\bfx) \d\bfx + f(\psi_{\obj}) + f(\psi_{\op}),
}
and, similarly, $\int g(\bfz) \d\bfz = \int_{\boZ^{\inx}} g(\bfz) \d\bfz + g(\phi)$ for any integrable function $g$ on $\boZ$.

Consider the hypothesis that a distribution $p$ on $\boX$ represents a true target. We define $\epsilon$ as the random variable on $\{0,1\}$ describing the fact that this hypothesis is true ($1$) or false ($0$) and introduce $\varphi$ as the virtual state occupied by a target whose existence is assumed by a false hypothesis. The conditional distribution $P(\cdot \given \epsilon)$ of the target therefore verifies
\eqns{
P(\cdot \given \epsilon = 0) = \delta_{\varphi} \spacedAnd P(\cdot \given \epsilon = 1) = p,
}
which induces a distribution $\bp$ on the extended state space $\bar\boX \defeq \boX \cup \{\varphi\}$ defined as
\eqnla{eq:expansionBarP}{
\bp(\bfx) & = P(\bfx \given \epsilon = 0) \Pr(\epsilon = 0) + P(\bfx \given \epsilon = 1) \Pr(\epsilon = 1) \nonumber \\
& = (1-w)\delta_{\varphi}(\bfx) + w p(\bfx)
}
for any $\bfx \in \bar\boX$, with $w = \Pr(\epsilon = 1)$ and with $\delta_{\varphi}$ the Dirac function at $\varphi$ such that $\delta_{\varphi}(\bfx)$ equals to $1$ if $\bfx = \varphi$ and $0$ otherwise (well defined since $\varphi$ is assumed to be an isolated point). This extensions is particularly useful for proving results since it combines the \emph{probability of existence} $w$ and the law $p$ into a single distribution $\bar{p}$ on $\bar\boX$.

\subsection{Sensor modelling}

The sensor is understood as a \emph{finite-resolution} sensor which can only generates observations in a finite partition $\Pi_t$ of $\boZ^{\inx}$ at time $t \in \bbT$. With each \emph{observation cell} in $\Pi_t$ is associated a unique index and the set of all these indices is denoted $Z'_t$. A set $\{A^z_t \st z \in Z_t\}$ of subsets of $\Pi_t$ with $Z_t \subseteq Z'_t$ is made available by the sensor at each time step $t \in \bbT$ and corresponds to the actual observation of the targets in the system. 
It is also helpful to consider another subset $A^{\phi}_t = \{\phi\}$ that does not correspond to the given data but which will be associated with undetected targets. The set of subsets corresponding to the observations at time $t$ is $\{ A^z_t \st z \in \bar{Z}_t \}$, where $\bar{Z}_t = Z_t \cup \{\phi\}$.

Although this sensor modelling is not the most usual in the target-tracking literature, many physical sensors, such as radars, sonars or cameras, are actually finite-resolution sensors. This approach will also motivate the introduction of a special form of likelihood in \cref{sec:alternativeSensorModelling}.

\subsection{Single-target modelling}
\label{ssec:individualModelling}

\subsubsection{Transition} 
\begin{figure}
\centering
\def\svgwidth{.75\columnwidth}
\scriptsize
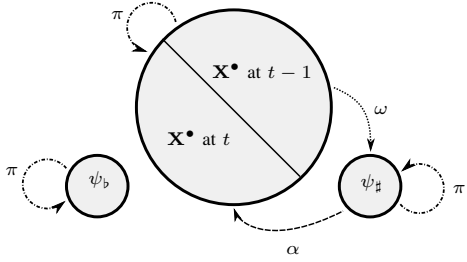
\caption{Relation between the subsets of $\boX$ between times $t-1$ and $t$ induced by the transition functions. The lines describe the transfer of probability mass from one subset to another, e.g.\ $q^{\alpha}_t$ transfers mass from $\psi_{\obj}$ to $\boX^{\inx}$ at time $t$ but to and from no other subsets.}
\label{fig:kernels}
\end{figure}

Three (sub-)transition functions $q^{\iota}_t$ from $\boX$ to itself indexed by $\iota \in \{\alpha,\pi,\omega\}$ are introduced in order to model the motion as well as the appearance and disappearance of targets between times $t-1$ and $t$. These functions are said to be sub-transitions because $q^{\iota}_t(\cdot \given \bfx)$ is not a probability distribution in general but is assumed to verify $\int q^{\iota}_t(\bfx' \given \bfx) \d\bfx' \leq 1$). Since there is no possible transition between the subset $\boX^{\inx} \cup \{\psi_{\obj}\}$ describing the targets and the point $\psi_{\op}$ describing the false-alarm generators, it holds for any $\iota \in \{\alpha,\pi,\omega\}$ that
\eqns{
q^{\iota}_t(\bfx' \given \psi_{\op}) = q^{\iota}_t(\psi_{\op} \given \bfx) = 0
}
for any $\bfx,\bfx' \neq \psi_{\op}$. These transition functions can then be characterised as follows:
\begin{enumerate}[label=\roman*),wide]
\item $q_t^{\alpha}$ models the appearance of a target, i.e.\ the transition from $\psi_{\obj}$ to $\boX^{\inx}$, so that $q_t^{\alpha}(\cdot \given \bfx) = \zero$ for any $\bfx \neq \psi_{\obj}$, where $\zero$ is the function equal to $0$ everywhere. It is assumed that a target cannot appear and disappear during one time step so $q_t^{\alpha}(\psi_{\obj} \given \psi_{\obj}) = 0$.
\item $q_t^{\pi}$ models targets' dynamics, i.e.\ transitions from $\boX^{\inx}$ to $\boX^{\inx}$ or from $\psi_{\obj}$ or $\psi_{\op}$ to themselves, so that for any $\bfx \in \boX^{\inx}$ it holds that
\eqns{
q_t^{\pi}(\psi_{\obj} \given \bfx) = 0 \spacedAnd q_t^{\pi}(\psi_{\obj} \given \psi_{\obj}) = q_t^{\pi}(\psi_{\op} \given \psi_{\op}) = 1.
}
\item $q_t^{\omega}$ models targets' disappearance, that is a transition from $\boX^{\inx}$ to $\psi_{\obj}$, so that for any $\bfx,\bfx' \in \boX^{\inx}$ it holds that
\eqns{
q_t^{\omega}(\bfx' \given \bfx) = 0 \spacedAnd q_t^{\omega}(\cdot \given \psi_{\obj}) = q_t^{\omega}(\cdot \given \psi_{\op}) = \zero.
}
\end{enumerate}
The transition functions $q_t^{\pi}$ and $q_t^{\omega}$ are additionally assumed to verify
\eqnl{eq:assumptionKernel}{
q_t^{\omega}(\psi_{\obj} \given \bfx) + \int q_t^{\pi}(\bfx' \given \bfx)\d\bfx' = 1
}
for any $\bfx \in \boX$, since a target with state $\bfx$ in $\boX^{\inx}$ can be either be propagated to $\boX^{\inx}$ with probability $p^{\pi}_t(\bfx) \defeq \int q_t^{\pi}(\bfx' \given \bfx) \d\bfx'$ or can disappear and be moved to $\psi_{\obj}$. It is convenient to model the appearance of targets with a transition function, e.g.\ if a rate of appearance is available then the probability of appearance $w^{\alpha}_t \defeq \int q^{\alpha}_t(\bfx' \given \psi_{\obj})\d\bfx'$ will depend on the duration of the considered time step, which is natural for a transition-related quantity. A graphical representation of the transfer of probability mass induced by the different transition functions is shown in \cref{fig:kernels}.

In order to use these transitions in the \gls{hisp}, they have to be extended to include the point $\varphi$ as well. For any $\iota \in \{\alpha,\pi,\omega\}$ the transition function $\bq^{\iota}_t$ from $\bar\boX$ to $\bar\boX$ is defined as
\eqnl{eq:extensionKernel}{
\bq^{\iota}_t(\bfx' \given \bfx) =
\begin{cases*}
q^{\iota}_t (\bfx' \given \bfx) & if $\bfx,\bfx' \in \boX$ \\
1-\int q^{\iota}_t (\bfx'' \given \bfx) \d\bfx'' & if $\bfx \in \boX$ and $\bfx' = \varphi$ \\
1 & if $\bfx = \varphi$ and $\bfx' = \varphi$.
\end{cases*}
}
This can be considered as the natural extension of a transition function from $\boX$ to $\boX$ to the corresponding extended spaces. The transition $\bq^{\alpha}_t$ has a \emph{multiplicity} $n^{\alpha}_t$, meaning that it is used exactly $n^{\alpha}_t$ times. It is assumed that the maximum number of appearing targets $n^{\alpha}_t$ is larger than the number of observations at time $t$ so that all observations might correspond to appearing targets. One possibility is to consider as many possibly appearing targets as there are resolution cells in the sensor. The transition $\bq^{\pi}_t$ and $\bq^{\omega}_t$ have unconstrained multiplicities so that their number simply adapt to the number of targets.\footnote{This can be made rigorous using the concept of partial information \cite{Houssineau2017,Houssineau2016_dataAssimilation} as discussed in \cite{Houssineau2015}}

\subsubsection{Likelihood} Let $\ell_t^{\d}$ be a likelihood from $\bar\boX$ to $\boZ$ describing the possible detection of targets and verifying
\eqns{
\ell_t^{\d}(\phi \given \psi_{\obj}) = \ell_t^{\d}(\phi \given \varphi) = 1.
}
The probability for a target with state $\bfx \in \boX^{\inx}$ to generate an observation is therefore $p^{\d}_t(\bfx) \defeq \int_{\boZ^{\inx}} \ell_t^{\d}(\bfz \given \bfx)\d\bfz$. The state $\psi_{\op}$ is the only point of $\bar\boX$ at which the likelihood $\ell_t^{\d}$ does not integrate to $1$ and we assume instead that $\ell_t^{\d}(\cdot \given \psi_{\op}) = \zero$ to ensure that false-alarm generator cannot be updated with it.

For each observation cell indexed by $z \in Z'_t$, a likelihood $\ell_t^z$ from $\bar\boX$ to $\boZ$ is introduced such that $\ell_t^z(\cdot \given \psi_{\op}) = w^{\op,z}_t \ind{A^z_t}$, with $\ind{A}$ the indicator of $A$, such that $\ell_t^z(\phi\given \psi_{\op}) = 1 - w^{\op,z}_t$, with $w^{\op,z}_t \in [0,1]$ the probability for a false alarm to be generated in the observation cell $z$, and such that $\ell_t^z(\cdot \given \bfx) = \zero$ for any $\bfx \neq \psi_{\op}$. Each likelihood in $\{\ell_t^z \st z \in Z'_t\}$ is then assumed to be used exactly once, i.e.\ each of them has multiplicity one. This model ensures that no more than one false alarm will appear in each resolution cell. If false alarms are independently identically distributed, a single likelihood function can be used to represent them all, with a multiplicity of $|Z'_t|$. The set $Z'_t \cup \{\d\}$ gathers all the possible indices for likelihoods.

Although this way of describing the dynamics and observation models is more sophisticated than usual, it allows for modelling a closed system where targets merely change state when appearing and disappearing and where all observations come from some objects, either a target or a false-alarm generator.

\subsection{Multi-target modelling}
\label{ssec:populationModelling}

\subsubsection[Standard MTT observation model]{Standard \gls{mtt} observation model} It is assumed that each subset $A^z_t$, $z \in Z_t$, is an observation which corresponds to one object, a target or a false-alarm generator. It is also assumed that objects cannot be related to several observations.

\subsubsection{Multi-target configuration} In general, the represented multi-target system might have some distinguishable parts, e.g.\ the previously-detected targets, and some indistinguishable parts, e.g.\ the never-detected targets. This aspect makes challenging the expression of the joint law of all these targets. We consider instead \emph{multi-target configurations} of the form $\bar\calP = \{(\bp_i,n_i)\}_{i \in \bbI}$ where $\bp_i$ is the law on $\bar\boX$ for the target(s) with index $i$ and $n_i$ is the associated multiplicity, i.e.\ there are $n_i$ targets with law $\bp_i$ in the considered configuration. It is assumed without loss of generality that $\bp_i \neq \bp_j$ for any $i,j \in \bbI$ such that $i \neq j$. If we assume that $\{\bp_i\}_{i \in \bbI}$ contains all the possible single-target laws, then varying $n_i$ gives all the possible configurations. It remains to define a suitable index set $\bbI$ to make sure that the associated multi-target configurations are suitable for the problem at hand. Note that we allow the multiplicity of the laws $\delta_{\psi_{\obj}}$, $\delta_{\psi_{\op}}$ and $\delta_{\varphi}$ to take values in the extended set of natural numbers $\bar\bbN = \bbN \cup \{\infty\}$.

Since $\varphi$ does not represent an actual target state, care must be taken when considering aspects related to distinguishability. Indeed, targets are indistinguishable if and only if their respective laws are equal when conditioned on the targets being in $\boX$, i.e.\ the targets represented by the random variables $X$ and $X'$ with respective laws $\bp$ and $\bp'$ on $\bar\boX$ are indistinguishable if and only if
\eqnl{eq:indistVphi}{
\bp(\bfx \given X \in \boX) = \bp'(\bfx \given X' \in \boX)
}
for all $\bfx \in \boX$.


\begin{example}
If $\bbI$ is a singleton then the multi-target configuration~$\bar\calP = \{(\bp,N)\}$ is equivalent to the symmetric multi-target law $P$ on $\bar\boX^N$ characterised for any $\bfx_1, \dots, \bfx_N \in \bar\boX$ by
\eqnsa{
P(\bfx_1, \dots, \bfx_{n_i}) & = \bp(\bfx_1) \dots \bp(\bfx_N) \\
& = \prod_{i=1}^N \Big((1-w)\delta_{\varphi}(\bfx_i) + w p(\bfx_i)\Big)
}
In this specific case, $P$ is symmetrical in its arguments and is related to the law of a multi-Bernoulli \gls{rfs}. Alternatively, if $\bbI =\{ 1,\dots,N \}$ and $n_i = 1$ for all $i \in \bbI$ then $\bar\calP$ is equivalent to the multi-target law $Q$ on $\bar\boX^N$ characterised for any $\bfx_1, \dots, \bfx_N \in \bar\boX$ by
\eqnsa{
Q(\bfx_1, \dots, \bfx_N) & = \bp_1(\bfx_1) \dots \bp_N(\bfx_N) \\
& = \prod_{i=1}^N \Big( (1-w_i)\delta_{\varphi}(\bfx_i) + w_i p_i(\bfx) \Big),
}
where each law $\bp_i$ has been expressed as in \cref{eq:extensionKernel}, which is similar to the law of a labelled multi-Bernoulli \gls{rfs} \cite{Vo2013}. The law $P$ does not distinguish any target whereas $Q$ distinguishes all targets. The objective is then to use other multi-target configurations which allow for partial distinguishability.
\end{example}

\begin{remark}
Even if labelled \gls{rfs} were generalised to permit several labels to be equal, the \gls{rfs} representation would not allow for a representation of indistinguishable targets almost surely at an isolated point state such as $\psi_{\obj}$ since a \gls{rfs} is a simple point process and therefore does not allow for multiplicity strictly greater than $1$ (both the state and the label would be equal in the case of indistinguishable targets at point~$\psi_{\obj}$). Handling non-simple, partially-labelled point processes is a complex task, and multi-target configurations alleviate much of this complexity by appropriately representing the hierarchical structure of partially-distinguishable multi-target system. A more general discussion about the different approaches for target tracking is provided in \cref{sec:relationWithOtherWorks}.
\end{remark}


\subsubsection{Observation path} With the standard \gls{mtt} observation model, the targets are made distinguishable as soon as they are detected since each observation corresponds to no more than one target. Also, one of the main characteristics of the induced single-target laws is the corresponding \emph{observation path} (the law is basically the Bayesian posterior law of one target state given the associated observations). For this reason, we consider the set $\bar\bbO_t$, defined as
\eqns{
\bar\bbO_t \defeq \bar{Z}_0 \times \dots \times \bar{Z}_t,
}
so that a sequence of observation $\bso_t \in \bar\bbO_t$ takes the form $\bso_t = (\phi,\ldots,\phi,z_{t_+},\ldots,z_{t_-},\phi,\ldots,\phi)$ with $t_+$ and $t_-$ the time of appearance and disappearance of the target, with $z_t \in \bar{Z}_t$ for any $t \in \{t_+,\dots,t_-\}$ and with $z_{t_+},z_{t_-} \neq \phi$. The empty observation path $(\phi,\dots,\phi) \in \bar\bbO_t$ is denoted $\bsphi_t$.

\subsubsection{Simplifying procedures} In order to reduce the number of terms in the considered multi-target configurations, two simplifying procedures are considered. They rely on \emph{mixing} two or more elements of a given multi-target configuration $\bar\calP$ on $\bar\boX$, say the elements with index $i$ and $j$ in the associated index set $\bbI$, by defining a new index $k$ based on $i$ and $j$ with $n_k = n_i + n_j$ and with the single-target law $\bp_k$ defined as
\eqns{
\bp_k(\bfx) = \dfrac{n_i \bp_i(\bfx) + n_j \bp_j(\bfx)}{n_i + n_j}.
}
for any $\bfx$ in $\bar\boX$. The index set after mixing $i$ and $j$ is then defined as $(\bbI \setminus \{i,j\}) \cup \{k\}$. The two simplifying procedures can then be formalised as
\begin{enumerate}[label=\textbf S.\arabic*]
\item \label{it:mixUndetected} Mixing of the appearing targets with the never-detected ones
\item \label{it:forgetDisappearedUndetected} Mixing with $\delta_{\psi_{\obj}}$ of the laws of never-detected targets that disappeared   
\end{enumerate}
Simplification~\ref{it:mixUndetected} fixes the number of groups of indistinguishable targets to one: the group of never-detected targets. The number of groups would otherwise grow by one every time step, with the appearance of additional targets. Simplification~\ref{it:mixUndetected} is justified by the fact that the distribution of appearing targets is often uninformative and constant in time, so that the difference between the appearing targets at time $t$ and the ones who appeared at time $t-1$ and who have not been detected at time $t$ can often be neglected. Similarly, Simplification~\ref{it:forgetDisappearedUndetected} is associated with the fact that never-detected disappeared targets are often irrelevant in \gls{mtt}. Their single-target law is the form $w\delta_{\psi_{\obj}} + (1-w)\delta_{\varphi}$ for some $w \in [0,1]$ that is often close to $0$, so that forcing the mixing with the single-target law $\delta_{\psi_{\obj}}$ only incurs a small information loss.

\subsubsection{Indexing of single-target laws} It is assumed that the only source of specific information at time ${t \in \bbT}$ lies in the observations made before time $t$. The targets' laws after prediction can then be indexed by the set $\bbI^{\obj}_{t|t-1}$ of triplets $(\obj,T,\bso)$ such that $T$ is either the empty set or a non-empty interval $[\cdot,t_-]$ of $\{0,\dots,t\}$ with unknown starting time and such that the observation time $\bso$ in $\bar\bbO_{t-1}$ satisfies $\bso_{t'} \neq \phi \implies t' \leq t_-$. An interval $[\cdot,t_-]$ is considered instead of a more standard given interval $\{t_+,\dots,t_-\}$ because newly appeared targets will be mixed with never-detected targets which implies that the time of appearance is forgotten. Although writing $t \in [\cdot,t_-]$ is equivalent to $t \leq t_-$, the expression $[\cdot,t_-]$ is preferred since it allows for using $\emptyset$ for targets that did not appear yet. The index set $\bbI^{\obj}_{t|t-1}$ can be interpreted as follows: predicted single-target laws are distinguished by their interval of presence in $\boX^{\inx}$ up to time $t$ and by their observation path in $\bar\bbO_{t-1}$. Simplification~\ref{it:forgetDisappearedUndetected} implies that the elements of the form $(\obj,[\cdot,t_-],\bsphi_t)$ with $t_- < t$ are not included in $\bbI^{\obj}_{t|t-1}$. The set $\bbI^{\obj}_t$ indexing single-target laws at time $t$ after the update is defined similarly but with the observation path in $\bar\bbO_t$, i.e.\ with the observation up to time $t$, rather than up to time $t-1$.

False-alarm generators are at $\psi_{\op}$ almost-surely and are assumed to give inconsistent observations, which corresponds to indices of the form $(\op,\emptyset,\bso)$ with $\bso$ containing a single non-empty observation. 
The set indexing the false-alarm generators at time $t$ after update is then
\eqns{
\bbI^{\op}_t = \bigcup_{0 \leq t' \leq t} \big\{ \big( \op,\emptyset,(\phi,\dots,\phi,z_{t'},\phi,\dots,\phi) \big) \st z_{t'} \in \bar{Z}_{t'} \big\}.
}

The sets indexing the multi-target configurations after prediction and update are then defined as
\eqns{
\bbI_{t|t-1} \defeq \bbI^{\obj}_{t|t-1} \cup \bbI^{\op}_{t-1} \spacedAnd \bbI_t \defeq \bbI^{\obj}_t \cup \bbI^{\op}_t.
}
The single-target laws, e.g.\ after the update, can then be indexed as follows:
\begin{enumerate}[label=\itshape\alph*\upshape)]
\item a target that is still present at time $t$ and that has only been detected during the current time step has index $(\obj,[\cdot,t],(\phi,\dots,\phi,z))$ for some $z \in Z_t$,
\item a target that has not appeared yet has index $(\obj,\emptyset,\bsphi_t)$ (which would also be an element of $\bbI^{\op}_t$ in the absence of the symbols $\obj$ and~$\op$), and
\item to each index $(\obj,T,\bso) \in \bbI^{\obj}_{t|t-1}$, corresponds $|\bar{Z}_t|$ indices in $\bbI^{\obj}_t$, that is one index $(\obj,T,\bso \times z)$ for each $z \in \bar{Z}_t$, with $\bso \times z$ denoting the concatenation of the sequence $\bso$ with the element $z$
\end{enumerate}
It is useful to introduce two more symbols to represent some specific parts of the multi-target system: the targets with an index $(\obj,T,\bso) \in \bbI^{\obj}_t$ that 
\begin{enumerate}
\item[($\m$)] have been previously detected, so that $\bso \neq \bsphi$, 
\item[($\u$)] are in the state space but have never been detected, so that $T \ni t$ and $\bso = \bsphi$,
\end{enumerate}
Because of Simplification~\ref{it:mixUndetected}, there is a single element in the set $\bbI^{\u}_t$, which we denote $\bsi^{\u}_t$ (similarly $\bbI^{\u}_{t|t-1} = \{\bsi^{\u}_{t|t-1}\}$). For the sake of compactness, single-target laws and multiplicities with index $\bsi^{\u}_t$ or $\bsi^{\u}_{t|t-1}$ will simply bear the superscript $\u$.

We will first consider multi-target configurations on $\bar\boX$ in the following sections, before reformulating some of the equations of the filter with multi-target configurations on $\boX$. In both cases, an element of the underlying collections will be referred to as an \emph{hypothesis}.

The following multi-target configuration represents the false-alarm generator, the yet-to-appear targets as well as some erroneous hypotheses:
\eqns{
\bar\calP^{\psi,\varphi}_t \defeq \{ (\delta_{\psi_{\op}}, n^{\bsi}_t) \}_{\bsi \in \bbI^{\op}_t} \cup \{ (\delta_{\psi_{\obj}}, n_{\psi_{\obj}}), (\delta_{\varphi}, n_{\varphi}) \},
}
where $n_{\psi_{\obj}} = \infty$ and $n^{\bsi}_t = n_{\psi_{\op}} = \infty$ when $\bsi = (\op,\emptyset,\bsphi_t)$. This means that there is an infinite number of targets that are not currently in the area of interest (the state space $\boX^{\inx}$) but which might enter it at a later time and there also are an infinite number of potential false-alarm generators (but the likelihoods corresponding to the actual generation of false alarms come in finite number). Setting $n_{\psi_{\obj}}$ and $n_{\psi_{\op}}$ equal to $\infty$ is a simplifying assumption as it implies that these cardinalities will never change during the scenario since subtracting any finite number to account for target appearance or the generation of false alarms will not affect them. The law $\delta_{\varphi}$ only serves as a representation of erroneous single-target laws, so that the value of $n_{\varphi}$ is irrelevant and its time dependency is omitted. The initial multi-target configuration $\bar\calP_0$ is defined as
\eqns{
\bar\calP_0 = \bar\calP_{0|-1} \defeq \{(\bp^{\u}_0, n^{\u}_0)\} \cup \bar\calP^{\psi,\varphi}_0
}
with $n^{\u}_0 = n^{\alpha}_0$, $\bp^{\u}_0(\bfx) = \bq^{\alpha}_0(\bfx \given \psi_{\obj})$ for any $\bfx \in \boX$.

\subsection{Prediction}
\label{ssec:timeFiltering}

The number of never-detected targets at time $t-1$ after the update is denoted $n^{\u}_{t-1}$. The multi-target configuration $\bar\calP_{t-1}$ after the update at time $t-1$ is assumed to have the following form:
\eqnl{eq:obsUpPrevTime}{
\bar\calP_{t-1} = \{ (\bp^{\bsi}_{t-1}, 1)\}_{\bsi \in \bbI^{\m}_{t-1}} \cup \{(\bp^{\u}_{t-1},n^{\u}_{t-1})\} \cup \bar\calP^{\psi,\varphi}_{t-1}.
}

Let $\sigma_{t|t-1} : \bbI_{t-1} \times \{\alpha,\pi,\omega\} \to \bbI_{t|t-1}$ be the one-to-one mapping giving an index in $\bbI_{t|t-1}$ for each pair of indices in $\bbI_{t-1} \times \{\alpha,\pi,\omega\}$ corresponding to an hypothesis at time $t-1$ and a transition function. This mapping can be defined as
\eqns{
\sigma_{t|t-1}: ( (a,T,\bso), \iota ) \mapsto
\begin{cases*}
(a,T \cup \{t\},\bso) & if $\iota \in \{\alpha,\pi\}$ \\
(a, T, \bso) & otherwise,
\end{cases*}
}
with $a$ begin either $\obj$ or $\op$, since the presence of an target is extended to time $t$ with the transitions $q^{\alpha}_t$ and $q^{\pi}_t$ only. The index $\sigma_{t|t-1}(\bsi,\iota)$ is simply the index of the hypothesis obtained when predicting the previous hypothesis $\bsi$ with the transition $\iota$. The prediction can now be expressed as in the following theorem, which proof is given in the appendix, and which is based on the following approximation.
\begin{enumerate}[label=\bfseries A.\arabic*,series=Approx]
\item \label{predApprox} The hypotheses formed by predicting each single-target law in $\bbI^{\obj}_{t-1}$ with the transitions $q^{\pi}_t$ and $q^{\omega}_t$ are independent.
\end{enumerate}

\begin{theorem}
\label{thm:timeFiltering}
Under \ref{predApprox}, the multi-target configuration $\bar\calP_{t|t-1}$ after prediction to time $t$ is characterised by
\eqns{
\bar\calP_{t|t-1} = \{(\bp^{\bsi}_{t|t-1},1)\}_{\bsi \in \bbI^{\m}_{t|t-1}} \cup \{(\bp^{\u}_{t|t-1}, n^{\u}_{t|t-1})\} \cup \bar\calP^{\psi,\varphi}_{t-1},
}
where $n^{\u}_t = n^{\u}_{t-1} + n^{\alpha}_t$ is the predicted number of never-detected targets and where the predicted single-target law $\bp^{\bsi}_{t|t-1}$ is defined on $\bar\boX$ for any index $\bsi \in \bbI_{t|t-1}$ with $(\bsk,\iota) = \sigma_{t|t-1}^{-1}(\bsi)$ as
\eqns{
\bp^{\bsi}_{t|t-1}(\bfx) = \dfrac{n^{\u}_{t-1} \int \bq^{\pi}_t(\bfx \given \bfx') \bp^{\u}_{t-1}(\bfx')\d\bfx' + n^{\alpha}_t \bq_t^{\alpha}(\bfx \given \psi_{\obj})}{n^{\u}_{t-1} + n^{\alpha}_t}
}
if $\bsk = \bsi^{\u}_{t|t-1}$ and as $\bp^{\bsi}_{t|t-1}(\bfx) = \int \bq^{\iota}_t(\bfx \given \bfx') \bp^{\bsk}_{t-1}(\bfx')\d\bfx'$ otherwise.
\end{theorem}

Note that $\sigma_{t|t-1}^{-1}(\bsi)$ simply recovers the indices of the previous hypothesis and transition that lead to $\bsi$. It appears that the predicted configuration $\bar\calP_{t|t-1}$ takes the same form as the posterior configuration $\bar\calP_{t-1}$ at time $t-1$ expressed in~\cref{eq:obsUpPrevTime}. Although the single-target laws of disappeared targets are not very informative, they can be useful in practice since the scalar $\int \bq^{\omega}_t(\psi_{\obj} \given \bfx) \bp^{\bsk}_{t-1}(\bfx) \d\bfx$ gives the credibility of the hypothesis that the target with index $\bsk$ at time $t-1$ disappeared between $t-1$ and $t$.

\subsection{Update}
\label{ssec:observationFiltering}

Let $\sigma_t : \bbI_{t|t-1}\times (Z'_t \cup \{\d\}) \times (Z'_t \cup \{\varphi\}) \to \bbI_t$ be the one-to-one mapping describing the connections of all possible combinations of prior single-target law, likelihood and observation with the indices in $\bbI_t$ be defined as
\eqns{
\sigma_t : ((a,T,\bso), s, z ) \mapsto (a,T,\bso \times z).
}
Note that, by construction, many combinations will result in a posterior with all its probability mass on $\varphi$, such as with applying the likelihood $\ell^{\d}_t$ to a target with state $\psi_{\op}$. 
To facilitate the statement of the update, let
\eqns{
L^{s,z}_t(\bfx) = \int \ind{A^z_t}(\bfz) \ell_t^s(\bfz \given \bfx) \d\bfz
}
be the probability that an object with state $\bfx \in \boX$ has generated an observation in the cell $z \in Z'_t$ under the likelihood $\ell^s_t$ with $s \in Z'_t \cup \{\d\}$. Similarly, we denote $p^{\bsk,s}_t(z)$ as the marginal likelihood of the observation $z$ when the prior $\bp^{\bsk}_{t|t-1}$, $\bsk \in \bbI_{t|t-1}$, is restricted to $\boX$, that is
\eqns{
p^{\bsk,s}_t(z) = \int_{\boX} L^{s,z}_t(\bfx) \bp^{\bsk}_{t|t-1}(\bfx) \d\bfx,
}
and $\bp^{\bsk,s}_t(z) = p^{\bsk,s}_t(z) + L^{s,z}_t(\varphi)\bp^{\bsk}_{t|t-1}(\varphi)$ as the extended marginal likelihood. Note that the value of $s$ is known once $\bsk$ is fixed, so that it is possible to omit it and write $\bp^{\bsk}_t(z)$ instead of $\bp^{\bsk,s}_t(z)$. The next theorem, proved in the appendix, relies on the following approximation.
\begin{enumerate}[Approx,resume]
\item \label{upApprox} The hypotheses formed by updating each single-target law in $\bbI_{t|t-1}$ by each observation in $\bar{Z}_t$ are independent.
\end{enumerate}

\begin{theorem}
\label{thm:HispObsFilter}
Under \ref{upApprox}, the multi-target configuration $\bar\calP_t$ after the update at time $t$ is characterised by
\eqnl{eq:thm:HispObsFilter:fullDistribution}{
\bar\calP_t = \{ (\bp^{\bsi}_t, 1)\}_{\bsi \in \bbI^{\m}_t} \cup \{(\bp^{\u}_t, n^{\u}_t)\} + \bar\calP^{\psi,\varphi}_t,
}
where the posterior single-target law $\bp^{\bsi}_t$ with index $\bsi \in \bbI^{\obj}_t$, denoting $(\bsk,s,z) \defeq \sigma_t^{-1}(\bsi)$, is characterised for any $\bfx \in \boX$ by
\eqns{
\bp^{\bsi}_t(\bfx) = w^{\bsi}_t \dfrac{L^{s,z}_t(\bfx) \bp^{\bsk}_{t|t-1}(\bfx) }{ p^{\bsk,s}_t(z) },
}
and by $\bp^{\bsi}_t(\varphi) = 1 - w^{\bsi}_t$, where $w^{\bsi}_t$ is defined as
\begin{subequations}
\label{eq:thm:HispObsFilter}
\eqnl{eq:perIndividual}{
w^{\bsi}_t = \dfrac{w_{\ex}^{\bsk, z} p^{\bsk,s}_t(z)}{\sum_{z' \in \bar{Z}_t} w_{\ex}^{\bsk , z'} \bp^{\bsk,s}_t(z')}
}
or, equivalently, as
\eqnl{eq:perObs}{
w^{\bsi}_t = \dfrac{w_{\ex}^{\bsk, z} p^{\bsk,s}_t(z)}{\sum_{\bsk' \in \bbI_{t|t-1}} w_{\ex}^{\bsk',z} \bp^{\bsk'}_t(z)},
}
\end{subequations}
and where $w_{\ex}^{\bsk,z}$ is the multi-target marginal likelihood corresponding to the update by the observations in ${Z_t \setminus \{z\}}$ of the targets with index in
\eqns{
\bbI^{\bsk}_{t|t-1} \defeq \big( \bbI^{\m}_{t|t-1} \setminus \{\bsk\} \big) \cup \{ \u, \op \},
}
and with multiplicity $n^{\u}_{t|t-1} - \ind{\u}(\bsk)$ for the never-detected targets.
\end{theorem}

Less formally, the marginal likelihood $w_{\ex}^{\bsk,z} \in [0,1]$ can be understood as the assessment of the compatibility between the predicted laws and the collection of observations at the current time excluding the/a target with index $\bsk$ and the observation~$z$. Note that we consider $Z_t \setminus \{z\}$ rather than $\bar{Z}_t \setminus \{z\}$ in the definition of $w_{\ex}^{\bsk,z}$ since the empty observation $\phi$ might be associated with an arbitrary number of targets, i.e.\ it is not because one target is not detected that other targets have to be detected.

An important feature of the \gls{hisp} can already be highlighted: an a posteriori probability of detection failure can be computed through \cref{eq:perIndividual} when $z = \phi$. Also, an a posteriori probability for an observation to be a false alarm is obtained when considering $\bsk \in \bbI^{\op}_{t|t-1}$. \Cref{thm:HispObsFilter} reveals the fact that the collection of marginalised single-target target laws $\{\bp^{\bsi}_t\}_{\bsi \in \bbI^{\obj}_t}$ can be seen as single-target filters in interaction, where the weights of the filters are $\{w^{\bsi}_t\}_{\bsi \in \bbI^{\obj}_t}$.

\section[Summary of the HISP filter recursion]{Summary of the \gls{hisp} recursion}
\label{sec:HISP:pdf}

The extension of the probability distributions and transition functions to $\bar\boX$ makes the proofs of the results easier, however, from an implementation viewpoint, an hypothesis is best described by a probability distribution $p$ on $\boX$ together with a scalar $w \in [0,1]$ describing the probability of the corresponding hypothesis to represent a target. This is the approach taken in \cite{Delande2017_AAS}, which is given here for completeness.

For any $\bsi \in \bbI^{\obj}_t$, the scalar $w^{\bsi}_t \in [0,1]$ and the probability distribution $p^{\bsi}_t$ on $\boX$ are defined by
\eqns{
w^{\bsi}_t \defeq \int_{\boX} \bp^{\bsi}_t(\bfx')\d\bfx' \spacedAnd p^{\bsi}_t(\bfx) = \bp^{\bsi}_t(\bfx) / w^{\bsi}_t,
}
for any $\bfx \in \boX$. The multi-target configuration can now be equivalently defined as $\calP_{t-1} = \{(p^{\bsi}_{t-1},w^{\bsi}_{t-1},n^{\bsi}_{t-1})\}_{\bsi \in \bbI_{t-1}}$.

\subsection{Approximations}

Approximations \ref{predApprox} and \ref{upApprox} are of the same nature: they allow hypotheses for \emph{co-existing}, i.e.\ for being part of a joint multi-target law, when it should not be the case. For instance, a target cannot simultaneously remain in $\boX^{\inx}$ \emph{and} disappear, and it cannot generate several observations under the standard \gls{mtt} assumptions. However, allowing for the co-existence of all hypotheses has for consequence the reduction of the number of multi-target configurations to one. There are two aspects in the \gls{hisp} that mitigate the effects of these approximations:
\begin{enumerate*}[label=\alph*)]
\item although all hypotheses are considered at once, each of them is allowed to be false (represented by $\varphi$), and
\item the propagation of the distinguishability enables an efficient track extraction to be devised, in which the standard assumptions of \gls{mtt} can be enforced.
\end{enumerate*}

\subsection{Prediction}

The hypothesis with index $\bsk = (\obj,[\cdot,t-1],\bso) \in \bbI^{\obj}_{t-1}$ can be propagated with the transition $q^{\pi}_t$ and yield the hypothesis $\bsi = (\obj,[\cdot,t],\bso) \in \bbI^{\obj}_{t|t-1}$ with weight and law
\eqnsa{
w^{\bsi}_{t|t-1} & = w^{\bsk}_{t-1} \int\!\!\!\int q^{\pi}_t(\bfx \given \bfx') p^{\bsk}_{t-1}(\bfx') \d \bfx'\d\bfx \\
p^{\bsi}_{t|t-1}(\bfx) & \propto \int q^{\pi}_t(\bfx \given \bfx') p^{\bsk}_{t-1}(\bfx') \d \bfx',
}
or it can disappear when propagated with the transition $q^{\omega}_t$, in which case the index in $\bbI_{t|t-1}$ of the corresponding hypothesis remains equal to $\bsk$ and the associated weight and law are
\eqnsa{
w^{\bsk}_{t|t-1} & =  w^{\bsk}_{t-1} \int q^{\omega}_t(\psi_{\obj} \given \bfx') p^{\bsk}_{t-1}(\bfx') \d\bfx' \\
p^{\bsk}_{t|t-1} & = \delta_{\psi_{\obj}}.
}
For never-detected targets, it holds that $n^{\u}_{t|t-1} = n^{\u}_{t-1} + n^{\alpha}_{t-1}$ and
\eqnsa{
p^{\u}_{t|t-1}(\bfx) & = \dfrac{n^{\u}_{t-1} \int q^{\pi}_t(\bfx \given \bfx') p^{\bsk}_{t-1}(\bfx') \d \bfx' + n^{\alpha}_{t-1} p^{\alpha}_t(\bfx)}{n^{\u}_{t-1} + n^{\alpha}_{t-1}} \\
w^{\u}_{t|t-1} & = \dfrac{n^{\u}_{t-1} w^{\bsk}_{t-1} + n^{\alpha}_{t-1} w^{\alpha}_t}{n^{\u}_{t-1} + n^{\alpha}_{t-1}},
}
with $p^{\alpha}_t = q^{\alpha}_t(\cdot \given \psi_{\obj})$ the distribution of the appearing targets.

\subsection{Update}

The hypothesis with index $\bsk = (\obj,[\cdot,t-1],\bso) \in \bbI^{\m}_{t-1}$ can be updated with the observation $z \in \bar{Z}_t$ via the likelihood $\ell^{\d}_t$ and yield the hypothesis $\bsi = (\obj,[\cdot,t],\bso\times z)$ with weight given by \cref{eq:thm:HispObsFilter} and law
\eqns{
p^{\bsi}_t(\bfx) \propto L^{\d,z}_t(\bfx) p^{\bsk}_{t|t-1}(\bfx).
}
The marginal likelihoods $p^{\bsk,\d}_t(z)$ and $\bp^{\bsk,\d}_t(z)$ appearing in \cref{eq:thm:HispObsFilter} can be expressed in this case as
\eqnsa{
p^{\bsk,\d}_t(z) & = w^{\bsk}_{t|t-1} \int L^{\d,z}_t(\bfx) p^{\bsk}_{t|t-1}(\bfx) \d\bfx \\
\bp^{\bsk,\d}_t(z) & = p^{\bsk,\d}_t(z) + \ind{\phi}(z)(1-w^{\bsk}_{t|t-1}).
}
Other hypotheses can be specified from \cref{thm:HispObsFilter} in the same way such as the ones corresponding to false alarms. For instance, the posterior probability for the observation $z \in Z_t$ to be a false alarm can be computed from \cref{eq:perObs} as
\eqns{
w^{\bsi}_t = \dfrac{w_{\ex}^{\op, z} w^{\op,z}_t}{\sum_{\bsk' \in \bbI_{t|t-1}} w_{\ex}^{\bsk',z} \bp^{\bsk'}_t(z)},
}
with $\bsi = (\op,\emptyset,\bsphi_{t-1}\times z)$. This quantity is not usually computed directly, but it will appear to be crucial in the proposed track-extraction procedure.

\subsection{Alternative sensor modelling}
\label{sec:alternativeSensorModelling}

In some situations, it is simpler to assume that observations can be any point $\bfz$ of the observation space $\boZ^{\inx}$ rather than a resolution cell. For instance, the shape of the resolution cells can be approximated by a Gaussian function of the form
\eqns{
f_{\bfz}(\bfz') = \exp\bigg( -\dfrac{1}{2} (\bfz' - \bfz)^T \Sigma^{-1} (\bfz' - \bfz) \bigg),
}
where $\Sigma$ approximates the extent of the corresponding resolution cell. In this case, the analogue of the function $L^{\d,z}_t$ takes the form
\eqns{
L_t(\bfz,\bfx) = \int f_{\bfz}(\bfz') \ell^{\d}_t(\bfz' \given \bfx) \d\bfz'.
}
If the noise in the propagation of the signal emitted by the target and received by the sensor is negligible when compared to the size of the resolution cells then it holds that $L(\bfz,\bfx) = f_{\bfz}(H(\bfx))$ where $H$ is the observation function, that is,
\eqns{
L_t(\bfz,\bfx) = \exp\bigg( -\dfrac{1}{2} (H(\bfx) - \bfz)^T \Sigma^{-1} (H(\bfx) - \bfz) \bigg).
}
The expression of $L_t$ is very close to the one of a standard likelihood function, except that there is no normalising constant ($L(\cdot,\bfx)$ has maximum $1$ for all $\bfx \in \boX$). This approach can also be justified through a direct modelling of uncertainty \cite{Houssineau2017,Houssineau2016_dataAssimilation} which has connections with Dempster-Shafer theory \cite{Dempster1968,Shafer1976}.

\section{Approximation of the multi-target marginal likelihood}
\label{sec:mainApproximation}

\Cref{thm:HispObsFilter} is based on the yet-to-be-defined multi-target marginal likelihood $w_{\ex}^{\bsk,z}$ which value is needed for all pairs $(\bsk,z)$ in $\bbI_{t|t-1} \times \bar{Z}_t$. However, the computation of these marginal likelihoods comes at the cost of a high complexity which, if performed exactly, would significantly reduce the efficiency of the proposed method. Instead, we formulate two possible approximations which are related to the ``sparsity'' of the scenario, either from the viewpoint of the targets or from the viewpoint of the observations: for all $\bsk \in I$ with $I \subseteq \bbI^{\m}_{t|t-1}$ and for all $z \in Z$ with $Z \subseteq Z_t$, we assume that
\begin{enumerate}[Approx,resume]
\item \label{approx1} $p^{\bsk, \d}_t(z)p^{\bsk', \d}_t(z) \approx 0$ for any $\bsk' \in I$ such that $\bsk \neq \bsk'$, or
\vspace{0.4em}
\item \label{approx2} $p^{\bsk, \d}_t(z) p^{\bsk, \d}_t(z') \approx 0$ for any $z' \in Z$ such that $z \neq z'$.
\end{enumerate}
Considering \ref{approx1} for a given $I$ and a given $Z$ is equivalent to assuming that two single-target laws with index in $I$ are unlikely to be successfully updated (in terms of marginal likelihood) with the same observation $z \in Z$.
Approximation~\ref{approx2} is the counterpart of \ref{approx1}, for which two observations in $Z$ are unlikely to successfully update the same single-target law $\bp^{\bsk}_t$ with $\bsk \in I$. These two approximations allow for factorising the expression of $w_{\ex}^{\bsk,z}$. The results will be given using \ref{approx1}, the analogue with \ref{approx2} follows a very similar path.

\begin{theorem}
\label{cor:HispFactWeights}
For any $\bsk \in \bbI_{t|t-1}$ and any $z \in \bar{Z}_t$, applying \ref{approx1} to the subsets $\bbI^{\bsk}_{t|t-1}$ and $Z_t \setminus \{z\}$, the scalar $w_{\ex}^{\bsk,z}$ can be factorised as follows
\eqns{
w_{\ex}^{\bsk, z} = C'_t(\bsk,z) \prod_{\bsk' \in\, \bbI^{\m}_{t|t-1} \setminus \{\bsk\}} \bigg[ \bp^{\bsk', \d}_t(\phi) + \sum_{z' \in Z_t \setminus \{z\}} \dfrac{ p^{\bsk', \d}_t(z')}{C^{\u,\op}_t(z')} \bigg]
}
where
\eqnl{lem:factPjoinT:eq:Czt}{
C^{\u,\op}_t (z) = \dfrac{p^{\u, \d}_t(z)}{\bp^{\u, \d}_t(\phi)} + \dfrac{p^{\op, z}_t(z)}{\bp^{\op, z}_t(\phi)},
}
and where
\eqnsml{
C'_t(\bsk,z) = \big[ \bp^{\u, \d}_t(\phi) \big]^{n^{\u}_t-\ind{\u}(\bsk)} \\
\times \bigg[ \prod_{z \in Z'_t\setminus Z_{\op}} \bp^{\op, z}_t(\phi) \bigg] \bigg[ \prod_{z \in Z_t\setminus\{z\}} C^{\u,\op}_t(z) \bigg]
}
with $Z_{\op}$ equal to $\{z\}$ when $\bsk = \op$ and $\emptyset$ otherwise.
\end{theorem}

\Cref{cor:HispFactWeights} is a direct consequence of \cref{lem:factoP} stated in the appendix. Note that the marginal likelihood $\bp^{\op, z}_t(\phi)$ used in the \lcnamecref{cor:HispFactWeights} is simply equal to $\ell_t^z(\phi \given \psi_{\op})$. An important property of the \gls{hisp} that appears in \cref{cor:HispFactWeights} is that all the terms $w_{\ex}^{\bsk,z}$ can be computed with a complexity of order $\calO(|\bbI_{t|t-1}||Z_t|)$, as demonstrated in \cref{alg:compWex}. The computation of all the terms $C'_t(\bsk,z)$ has a lower complexity ($\calO(|Z_t|)$) and is not detailed.

\begin{algorithm}
\caption{Computation of $w^{\bsk,z}_{\ex}$ for all $(\bsk,z) \in \bbI_{t|t-1} \times \bar{Z}_t$}
\label{alg:compWex}
\begin{algorithmic}
    \FOR{$\bsk \in \bbI_{t|t-1}$}
        \FOR{$z \in Z_t$}
            \STATE $u_{\bsk,z} \leftarrow p^{\bsk, \d}_t(z) / C^{\u,\op}_t(z)$
        \ENDFOR
        \STATE $w_{\bsk,\phi} \leftarrow \bp^{\bsk, \d}_t(\phi)$
        \FOR{$z \in Z_t$}
            \STATE $w_{\bsk,\phi} \leftarrow w_{\bsk,\phi} + u_{\bsk,z}$
        \ENDFOR
        \FOR{$z \in Z_t$}
            \STATE $w_{\bsk,z} \leftarrow w_{\bsk,\phi} - u_{\bsk,z}$
        \ENDFOR
    \ENDFOR
    \FOR{$z \in \bar{Z}_t$}
        \STATE $W_z \leftarrow 1$
        \FOR{$\bsk \in \bbI_{t|t-1}$}
            \STATE $W_z \leftarrow W_z w_{\bsk,z}$
        \ENDFOR
    \ENDFOR
    \FOR{$\bsk \in \bbI_{t|t-1}$}
        \FOR{$z \in \bar{Z}_t$}
            \STATE $w^{\bsk,z}_{\ex} \leftarrow C'_t(\bsk,z) W_z / w_{\bsk,z}$
        \ENDFOR
    \ENDFOR
\end{algorithmic}
\end{algorithm}

\section{Relation with other works}
\label{sec:relationWithOtherWorks}

In this section, the relation between the proposed approach and other \gls{mtt} techniques is discussed.

\subsubsection*{\gls{lmb} filter} Distinguishing targets has been made possible with point processes by using labelling \cite{Vo2013,Vo2014}. The \gls{lmb} filter \cite{Reuter2014} follows as an approximation and is close in principle to the proposed approach. However, point processes have been built on the assumption that the targets are indistinguishable \cite[p.\ 124]{Daley2003} and labelling is usually meant to represent characteristics of the target that do not evolve in time instead of representing a target identity. These facts do not make the use of labels for distinguishing targets straightforward and specific techniques have to used to prevent the natural symmetrisation of point-process laws. The objective with the proposed framework is to build on a natural representation of partially-distinguishable multi-target systems \cite{Houssineau2016} which is based on a constructive approach leading in the independent case to a point process on the space of probability measures, and reducing to a multi-target configuration in the specific case of the \gls{hisp}. From a practical point of view, the labelled multi-Bernoulli filter computational complexity can be, in a worst-case scenario, as high as with non-approximated techniques, whereas the \gls{hisp}'s complexity is linear. 

\subsubsection*{Poisson multi-Bernoulli filter} The idea of separating the never-detected targets from the detected ones has been proposed in \cite{Williams2012,Williams2015} where detected and never-detected targets are respectively represented by a Poisson point process and by a (mixture of) multi-Bernoulli point process(es). The use of a Poisson distribution for appearing targets offers a practical advantage when no upper bound is known for the associated cardinality, as opposed to the case of a finite-resolution sensor considered here where the number of appearing targets cannot exceed the number of resolution cells. Since it is based on unlabelled point processes, the Poisson multi-Bernoulli approach does not allow for a principled track extraction or for any post-processing requiring targets to be distinguishable, such as classification, unlike the proposed \gls{mtt} algorithm \cite{Pailhas2016}. To reduce the mixture of multi-Bernoulli point processes obtained after the update, \cite{Williams2015} proposes to select the multi-Bernoulli distribution that minimises the Kullback-Leibler divergence with the mixture. This approach is well suited to unlabelled point processes where there is no question of distinguishability.

Many other approaches exist in the target-tracking literature, e.g.\ \cite{Garcia2016} considers random finite sets of trajectories.

\section{Implementation}
\label{sec:implementation}

Although the complexity is linear in the number of considered single-target laws and in the number of observations, specific approximations have to be used in practice to limit the computational cost and the number of propagated laws while ensuring that track extraction can be efficiently applied. 

Let $\bar\bbJ_{t-1} \subseteq \bbI^{\obj}_{t-1}$ contain the indices that have been retained up to time $t-1$ and denote $\bbJ_t \subseteq  \bbI^{\obj}_t$ the set of indices obtained at time $t$ after applying prediction and update to $\bar\bbJ_{t-1}$ (in the same way $\bbI_t$ is deduced from $\bbI_{t-1}$). Allowing different hypotheses to share the same single-target law, we denote $\bar\bbL_{t-1}$ the partition of $\bar\bbJ_{t-1}$ characterising this aspect at time $t-1$, that is, the index $\bsk \in \bar\bbL_{t-1}$ of a single-target law $p^{\bsk}_{t-1}$ is actually defined as the set containing the indices of the hypotheses described by the law $p^{\bsk}_{t-1}$. The objective with this approach is to allow for the merging of single-target laws while keeping distinct the other characteristics of the involved hypotheses. In this way, a large number of hypotheses can be propagated with a reasonable computational cost (which is mostly determined by the number of single-target laws to predict and update).

Let $\bbL_t$ be the propagated version of the partition $\bar\bbL_{t-1}$. Once again, an index $\bsk \in \bbL_t$ corresponds to a single-target law toward which several indices in $\bbJ_t$ can point. The weight $w^{\bsk}_t$ associated with the single-target law $\bsk \in \bbL_t$ is the sum of the weights of the hypotheses which relies on it, i.e.\ $w^{\bsk}_t = \sum_{\bsi \in \bsk} w^{\bsi}_t$ (which can be larger than $1$, as opposed to hypotheses' weight).

\subsubsection*{Pruning} Some hypotheses' weight will be very close to $0$ so that their probability of existence is low enough to discard them. The actual index set is then a subset $\tilde\bbJ_t$ of $\bbJ_t$. This \emph{pruning} procedure is characterised by
\begin{enumerate}[Approx,resume]
\item \label{it:pruning} The set $\tilde\bbJ_t$ is the subset of $\bbJ_t$ containing indices $\bsi$ such that $w^{\bsi}_t > \tau$ only.
\end{enumerate}
The procedure affects the indices of single-target laws, however we also denote by $\bbL_t$ the corresponding partition of $\tilde\bbJ_t$. Single-target laws will be automatically discarded when their index becomes empty (as a subset of $\tilde\bbJ_t$).

\subsubsection*{Merging} Some of the single-target laws will, in practice, be too close to each other to justify propagating them separately, a partition $\bar\bbL_t$ of $\tilde\bbJ_t$ can be introduced to group the laws that are alike. This is the \emph{merging} procedure characterised in the Gaussian case as follows (denoting $\calN(\bfm,\bsV)$ the Gaussian distribution with mean $\bfm$ and variance $\bsV$ and assuming $p^{\bsk}_t = \calN(\bfm^{\bsk}_t, \bsV^{\bsk}_t)$ for any $\bsk \in \bbL_t$).
\begin{enumerate}[Approx,resume]
\item \label{it:merging} The partition $\bar\bbL_t$ of $\tilde\bbJ_t$ is introduced recursively as:
\begin{enumerate}[label=(\roman*)]
\item Define $\bar\bbL_t$ on $K = \emptyset$ as the empty partition
\item \label{it:startRecMerging} Find the index $\bsk = \argmax_{\bsj \in \bbL_t \setminus K} w^{\bsj}_t$ corresponding to the single-target law with highest weight among the ones that have not already been merged and define $K'$ as the set containing any index $\bsj \in \bbL_t$ such that the Mahalanobis distance \cite{Mahalanobis1936} between $\calN(\bfm^{\bsk}_t, \bsV^{\bsk}_t + \bsV^{\bsj}_t)$ and $\bfm^{\bsj}_t$ is strictly less than $\tau'$ (recalling that each index $\bsj \in K'$ is itself a subset of $\tilde\bbJ_t$)
\item Let $\bsk'$ be the union of the cells in $K'$ and let $p^{\bsk'}_t$ be characterised by its mean and variance as
\eqnsa{
\bfm^{\bsk'}_t & = \dfrac{1}{\sum_{\bsk \in K'} w^{\bsk}_t} \sum_{\bsk \in K'} w^{\bsk}_t \bfm^{\bsk}_t \\
\bsV^{\bsk'}_t & = \dfrac{1}{\sum_{\bsk \in K'} w^{\bsk}_t} \\
& \quad \times \sum_{\bsk \in K'} w^{\bsk}_t \big( \bsV^{\bsk}_t + ( \bfm^{\bsk'}_t - \bfm^{\bsk}_t)(\bfm^{\bsk'}_t - \bfm^{\bsk}_t)^T \big)
}
\item Extend $\bar\bbL_t$ to $K \cup \bsk'$ by letting $\bsk'$ be a cell of the partition
\item Redefine $K$ as $K \cup \bsk'$ and return to step \ref{it:startRecMerging} until $K = \bbL_t$
\end{enumerate}
\end{enumerate}
The two indexed families of interest are then $\{ w^{\bsi}_t \}_{\bsi \in \tilde\bbJ_t}$ and $\{ p^{\bsk}_t \}_{\bsk \in \bar\bbL_t}$. These approximations are usual for handling Gaussian mixtures \cite{Salmond1990} but can can be applied here for any implementation of the filter by adapting the considered distance, e.g.\ the Hellinger distance \cite{Hellinger1909}.

Although the reduction of the number of hypotheses does not induce a computational gain as significant as the reduction of the number of single-target laws, there is still some interest in mixing hypotheses that are very similar, especially for long and/or complex scenarios. To decide when hypotheses are similar enough to be mixed, we consider a time window $T = \{t-l,\dots,t\}$ at time $t$ for some lag $l$ and require the corresponding observation paths to be the same over the window $T$. If there is a subset $I$ of hypotheses' indices verifying this assumption and if it holds that $w^I_t = \sum_{\bsi \in I} w^{\bsi}_t \leq 1$ then these hypotheses can be mixed: the observation path of the hypothesis with highest weight can be kept, e.g.\ for display purposes, and the weight of the new hypothesis is $w^I_t$. If the laws associated with each hypothesis being mixed are distinct then the resulting hypothesis is associated with a weighted mixture of the corresponding laws (although this does not typically happen when $l > 1$). We denote by $\bar\bbJ_t$ the set of hypotheses at time $t$ resulting from this mixing procedure.

\subsubsection*{Track extraction} As mentioned in previous sections, a posterior probability for an observation to be a false alarm is computed and the result is stored as an hypothesis for the purpose of track extraction. The hypotheses corresponding to disappeared targets are kept for the same reasons. The track extraction also operates on the time window~$T$, so that these hypotheses can be discarded once the time of their last observation falls out of this window. Finally, in order to perform track extraction, one can solve the following optimisation problem:
\eqnl{eq:optimise}{
\argmax_{I \subseteq \bar\bbJ_t} \; \prod_{\bsi \in I} w^{\bsi}_t
}
subject to:
\begin{enumerate}
\item \label{it:optimise1} the union of all observation paths over the time window $T$ must contain all the observations over this window,
\item \label{it:optimise2} the observations paths in $I$ must be \emph{compatible}: each observation cannot be used more than once.
\end{enumerate}
The solution to this problem is the same as the one for
\eqns{
\argmax_{I \subseteq \bar\bbJ_t} \; \sum_{\bsi \in I} \log w^{\bsi}_t
}
with the same constraints, since all $w^{\bsi}_t$ are strictly positive. The latter problem can however be solved by linear programming. Constraint~\ref{it:optimise1} justifies the fact that false alarms and disappeared targets are kept as hypotheses; the solution to \cref{eq:optimise} would not be meaningful otherwise. Constraint~\ref{it:optimise2} ensures that the assumptions of standard \gls{mtt} are satisfied. The only parameter for track extraction is the size of the time window $T$ which is practically appealing since it is easy to interpret and tune.

\begin{remark} To ensure that previously displayed tracks do not disappear when they have not been detected over the time window $T$, the corresponding observations can be kept even when their time is prior to $T$, i.e.\ observations corresponding to confirmed tracks are held longer in order to improve the results with a limited impact on the computational time.
\end{remark}

The track extraction procedure proposed in this section is only one among many possible. The fact that the \gls{hisp} provides distinct hypotheses enables the introduction of tailored extraction methods depending on the application and computational resources at hand. The procedure proposed in this section is considered since it is one of the simplest that utilises the structure of the filter as opposed, for instance, to selecting single-target laws based on their weight.

\section{Simulations}
\label{sec:simulations}

The performance of the \gls{hisp} is compared against the \gls{phd} and \gls{cphd} filters \cite{Mahler2007_CPHD} as well as the \gls{lmb} filter. Note that because of its hierarchical nature, the \gls{hisp} can be easily implemented using any Bayesian filtering technique for each single-target law, e.g.\ \gls{smc} as in \cite{Houssineau2015_SMC} or \gls{kf}.

We consider a sensor placed at the centre of the $2$-D Cartesian plane that delivers range and bearing observations every $\Delta = 4\si{s}$ during $200\si{s}$, i.e.\ the time index set is $\bbT \defeq \{0,\dots,50\}$ with the actual time being $4t$ for any $t \in \bbT$. The size of the resolution cells of this sensor is $1\si{\degree} \times 15\si{m}$. Considering small fixed random error and bias error, the standard deviation of the observations is $\sigma_r = 6.2\si{m}$ for the range and $\sigma_{\theta} = 4.5\si{mrad}$ for the bearing, for a \gls{snr} of $3\si{dB}$ and $\sigma_r = 4.87\si{m}$ and $\sigma_{\theta} = 3.5\si{mrad}$ for a \gls{snr} of $5\si{dB}$. The range $r$ is in $[50\si{m},500\si{m}]$ and the bearing $\theta$ is in $(-\pi,\pi]$. For the comparison with the (C)\gls{phd} and \gls{lmb} filters to be possible, point observations are generated according to the standard observation model with the standard deviations given above, instead of using the resolution cells. The alternative sensor modelling of \cref{sec:alternativeSensorModelling} is thus used for the \gls{hisp}.

The scenario comprises $5$ targets which motion is driven by a nearly-constant velocity model so that the random variable $X_t$ representing the state of a target in $\boX^{\inx}$ at time $t$ given its state $\bfx_{t-1}$ at the previous time verifies $X_t \sim \calN(\bsF \bfx_{t-1},\bsP)$ with
\eqns{
\bsF =
\begin{bmatrix}
1 & 0 & \Delta & 0 \\
0 & 1 & 0 & \Delta \\
0 & 0 & 1 & 0 \\
0 & 0 & 0 & 1
\end{bmatrix}
\text{, }
\bsP = \sigma^2
\begin{bmatrix}
\sfrac{\Delta^3}{3} & 0 & \sfrac{\Delta^2}{2} & 0 \\
0 & \sfrac{\Delta^3}{3} & 0 & \sfrac{\Delta^2}{2} \\
\sfrac{\Delta^2}{2} & 0 & \Delta & 0 \\
0 & \sfrac{\Delta^2}{2} & 0 & \Delta
\end{bmatrix}
}
with $\sigma^2 = 0.05\si{m^2/s^{4}}$. All targets are present at the beginning of the scenario and Targets~2 to 5 never spontaneously disappear whereas Target~1 disappear at $160\si{s}$ in Case~1 below and does not disappear in Cases~2 and 3. The scenario is depicted in \cref{fig:scenario}. Note that Targets~2 and 3 are crossing around $t = 120\si{s}$.

We consider a \gls{kf} implementation of the \gls{hisp} based on \ref{approx1}, \ref{it:pruning} and \ref{it:merging} and similarly for the \gls{lmb} filter. In this implementation, the detected and never-detected hypotheses are updated through \cref{eq:perObs} and \cref{eq:perIndividual} respectively. For the (C)\gls{phd} filter, a Gaussian mixture implementation\footnote{The code of the \gls{cphd} and \gls{lmb} filters was downloaded from Ba Tuong Vo's page (\url{http://ba-tuong.vo-au.com/rfs_tracking_toolbox_beta.7z})} \cite{Vo2006,Vo2007} with a confirmation threshold $\tau_{\c} = 0.9$ is considered, i.e.\ all the terms in the Gaussian mixture with a weight exceeding $\tau_{\c}$ are declared as confirmed tracks. The non-linearity of the observation model is dealt with by an extended Kalman filter. To reduce the computational cost, pruning (with parameter $\tau = 10^{-5}$) and merging (with parameter $\tau' = 4$) are carried out on the collection of posterior single-target laws or on the Gaussian mixture, depending on the filter. The probability for a target with state $\bfx \in \boX^{\inx}$ at time $t-1$ of remaining within $\boX^{\inx}$ at time $t$ is set to $p_{\pi} \defeq p^{\pi}_t(\bfx) = 0.999$.

In the considered scenarios, the mean number of appearing targets $m_{\alpha}$ is set to $10^{-2}$ per time step. Targets are assumed to be detected upon appearance in the \gls{phd} and \glspl{hisp}, the corresponding distribution is induced by the observation as in \cite{Houssineau2010} with the velocity components being initialised as Gaussian with mean $0$ and standard deviation $3.5\si{m/s}$; the associated weight is $w^{\alpha}_t$ which is defined as $m_{\alpha}/|Z'_t|$ for any $t \in \bbT$. The \gls{cphd} and \gls{lmb} filters are initialised with a Gaussian mixture whose terms are centred around the location of appearance of the targets, with a standard deviation of $50\si{m}$ on the position, and with the other parameters being the same as for the other filters. The average number of false alarms per time step is denoted~$n_{\op}$. The probability of detection is constant across the state space and through time, so that it is denoted $p_{\d} \defeq p^{\d}_t(\bfx)$ for any $\bfx \in \boX^{\inx}$ and any time $t$. From the given characteristics of the sensor and for a given value of $p_{\d}$, we deduce the probability for a single observation cell to produce a false alarm and we denote it $w_{\op}$. The approximate value of $n_{\op}$ can then be deduced directly from the number of observations cells. We proceed to the performance assessment on three different scenarios.

\subsubsection*{Case 1: High probability of detection $(5\si{dB})$}

We set $p_{\d} = 0.995$ so that $w_{\op} = 7.67\times 10^{-3}$ and $n_{\op} \approx 83$. The HISP track-extraction window is set to a length of $5$. The OSPA distance \cite{Schuhmacher2008} depicted in \cref{fig:case3} is based on a cutoff of $100$ and a $2$-norm and is averaged over $100$ \gls{mc} runs. Even though the estimation problem is not challenging with these parameters, there is a noticeable difference of performance between the two types of filters. This is mainly caused by the additional weighting term $w^{\bsk,z}_{\ex}$ of the \gls{hisp} which allows for a better discrimination between likely and unlikely hypotheses and which reduces the effects of association uncertainty on the overall performance. The \gls{cphd} filter takes the longest time to react to the disappearance of Target~1. The performance of the \gls{phd} filter is reduced when Targets~2 and 3 cross whereas the performance of the other filters is not affected.

\subsubsection*{Case 2: Low probability of detection $(3\si{dB})$}

We set $p_{\d} = 0.5$ so that $w_{\op} = 1.34\times 10^{-3}$ and $n_{\op} \approx 15$. The HISP track-extraction window is set to a length of $6$ since there might be many consecutive detection failures. The average OSPA distance is shown in \cref{fig:case1}. The OSPA distance for the \gls{hisp} is the lowest at all time steps. Due to the combination of a low probability of detection and the uncertainty on the association, the OSPA distance for the \gls{hisp} increases when Targets~$3$ and $4$ cross. The performance of the \gls{hisp} in this case is mainly explained by the fact that it computes an a posteriori probability of detection, so that the prior probability, $p_{\d} = 0.5$ here, has a lower impact on the final result when compared to the other methods.

\subsubsection*{Case 3: High probability of false alarms $(3\si{dB})$}

\begin{figure}
    \centering
    \begin{subfigure}[b]{0.806\columnwidth}
        \includegraphics[trim = 100pt 268pt 120pt 285pt,clip,width=\textwidth]{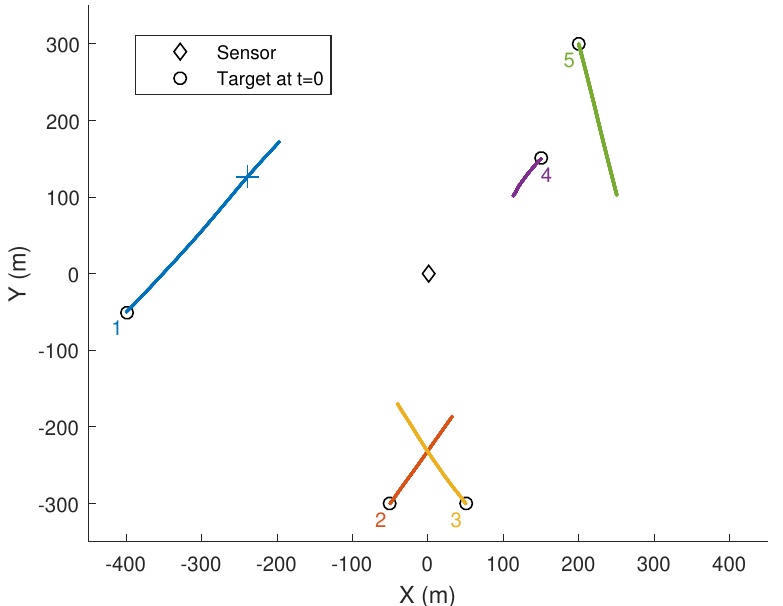}
        \caption{A realisation of the target trajectories (blue cross: location of Target~1 when it disappears in Case~1.)}
        \label{fig:scenario}
    \end{subfigure}\\
    \vspace{5pt}
    \begin{subfigure}[b]{0.804\columnwidth}
        \includegraphics[trim = 60pt 240pt 70pt 240pt,clip,width=\textwidth]{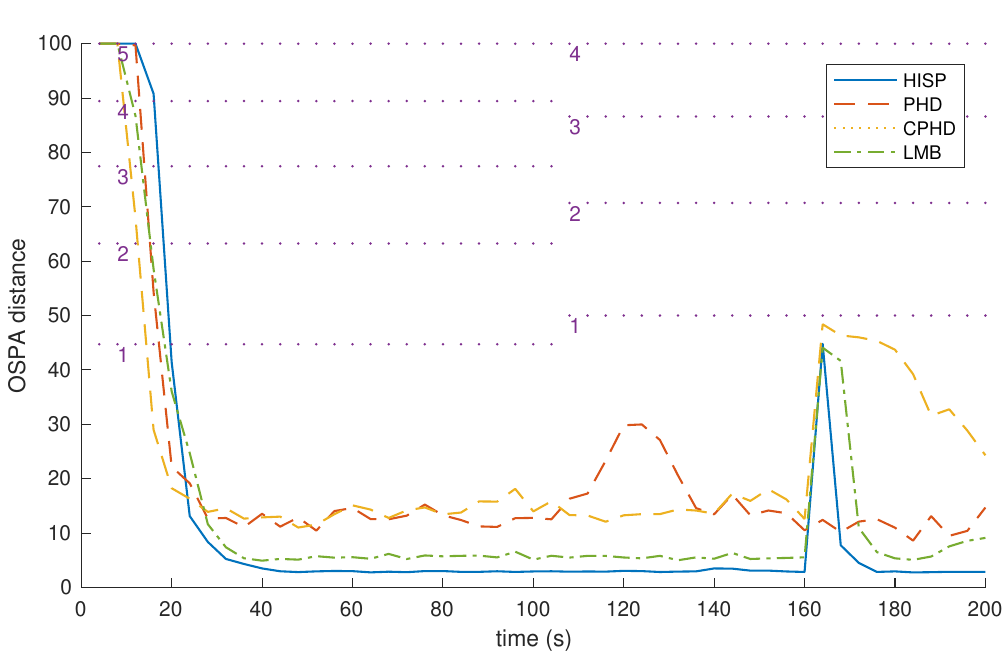}
        \caption{Case~1: $p_{\d} = 0.995$ and $n_{\op} \approx 83$.}
        \label{fig:case3}
    \end{subfigure}\\ 
    \vspace{5pt}
    \begin{subfigure}[b]{0.804\columnwidth}
        \includegraphics[trim = 60pt 240pt 70pt 240pt,clip,width=\textwidth]{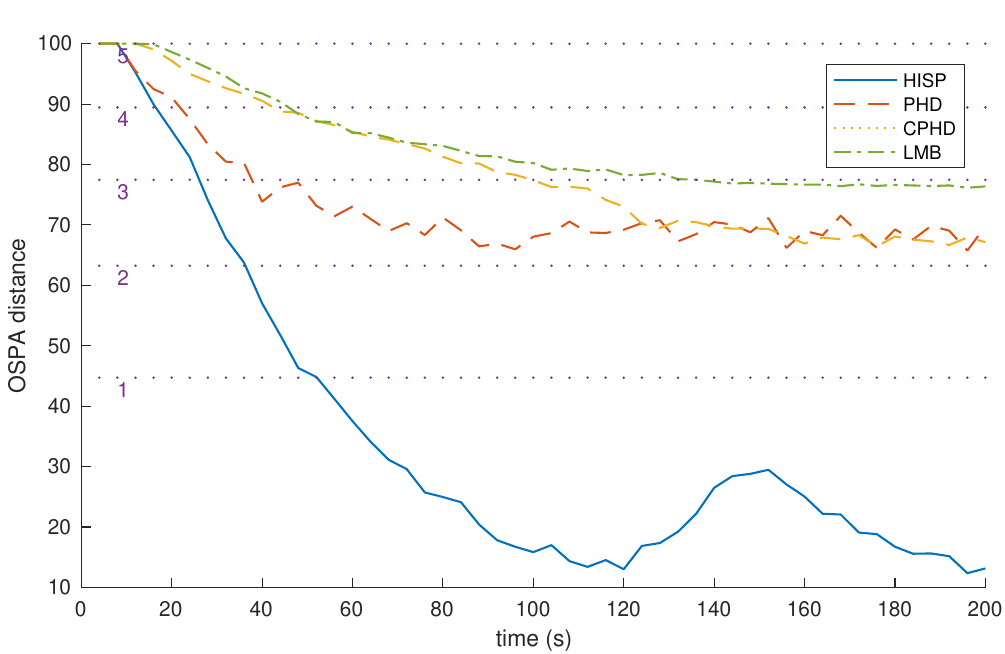}
        \caption{Case~2: $p_{\d} = 0.5$ and $n_{\op} \approx 15$.}
        \label{fig:case1}
    \end{subfigure}\\
    \vspace{5pt}
    \begin{subfigure}[b]{0.804\columnwidth}
        \includegraphics[trim = 60pt 240pt 70pt 240pt,clip,width=\textwidth]{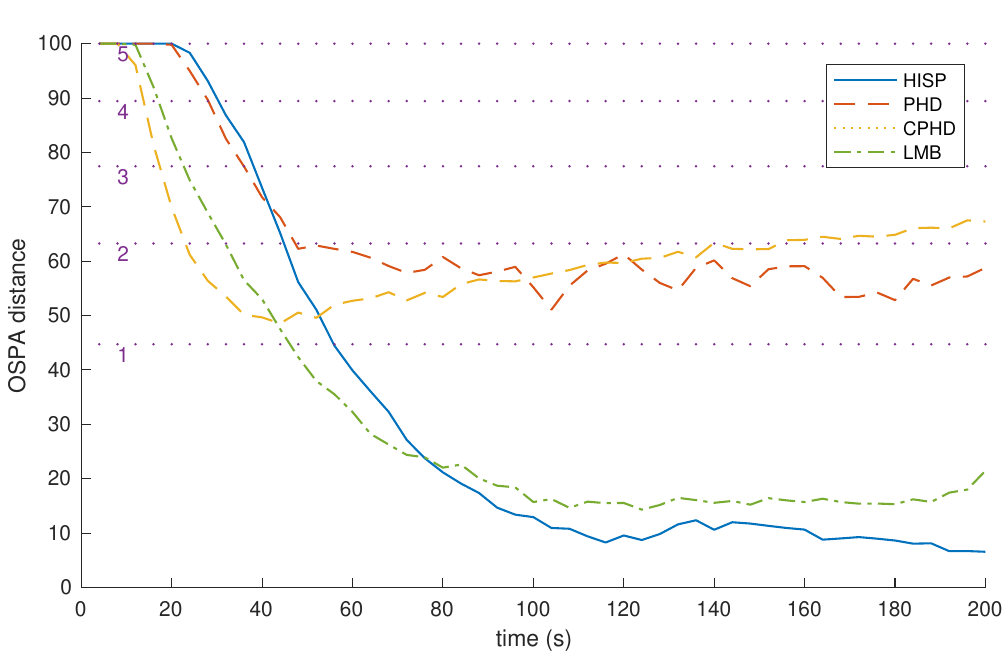}
        \caption{Case~3: $p_{\d} = 0.8$ and $n_{\op} \approx 167$.}
        \label{fig:case2}
    \end{subfigure}
    \caption{OSPA distance in Cases~1-3 (b-d) on the scenario (a) over $100$ \gls{mc} runs. (\gls{hisp}: solid line -- \gls{phd}: dashed line -- the dotted line numbered $n$ represents the OSPA for a cardinality-only error of $n$)}
\end{figure}

In this case, we set $p_{\d} = 0.8$ so that $w_{\op} = 1.54\times 10^{-2}$ and $n_{\op} \approx 167$. The HISP track-extraction window is set to a length of $3$ in order to alleviate the computational cost. The average OSPA distance is shown in \cref{fig:case2}. The \gls{phd} filter, which is known to be robust to high numbers of false alarms, behaves slightly better than in Case~2. The \gls{cphd} and \gls{lmb} filters react significantly faster to target appearance than the \gls{hisp} but this might be related to the more informative birth process they use. The \gls{cphd} filter is more prone to the creation of false tracks that impede its performance in the longer run.

\subsubsection*{Parameter sensitivity}

\begin{figure}
\centering
\includegraphics[trim = 100pt 265pt 120pt 285pt,clip,width=.75\columnwidth]{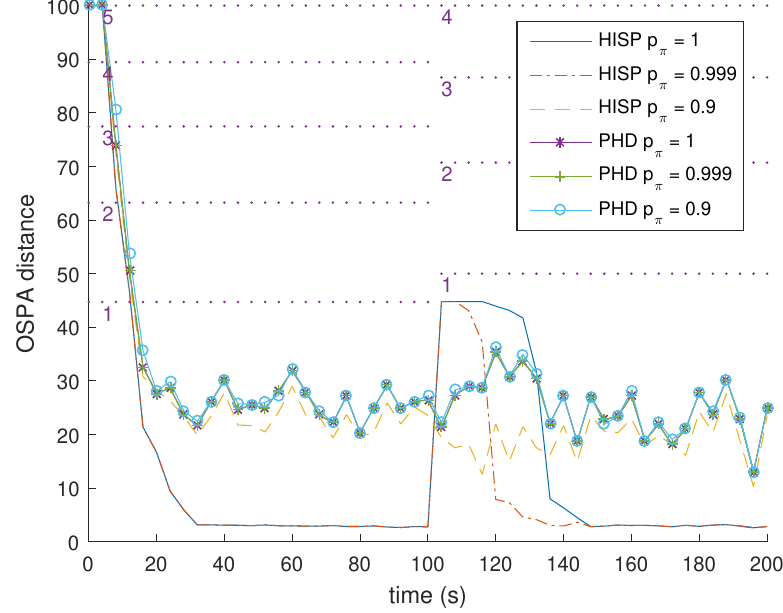}
\caption{OSPA distance versus time for different values of $p_{\pi}$, averaged over $50$ \gls{mc} runs.}
\label{fig:sensitivityPtPi}
\end{figure}

The \gls{hisp} displays a high sensitivity to some parameters when compared to the \gls{phd} filter. In particular, and as shown in \cref{fig:sensitivityPtPi}, the value of the probability $p_{\pi}$ has an important impact on the behaviour of the filter: $p_{\pi} = 1$ implies that if an hypothesis is almost-surely correct then it will be displayed at all following time steps, alternatively, if $p_{\pi} \leq p_{\d}$ then hypotheses stop to be considered as tracks as soon as a detection failure happens. Conversely, the behaviour of the \gls{phd} filter is nearly independent of $p_{\pi}$, so that this filter does not actually allow for taking the knowledge about persistence of targets into account. The scenario considered in \cref{fig:sensitivityPtPi} is a slightly modified version of the one considered above for Case~1, with a probability of detection $p_{\d} = 0.9$, with $n_{\op} = 10$ and with Target~1 disappearing at $100\si{s}$ rather than at $160\si{s}$.

\subsubsection*{Computational time}

\begin{figure}
\centering
\includegraphics[trim = 100pt 330pt 120pt 355pt,clip,width=.805\columnwidth]{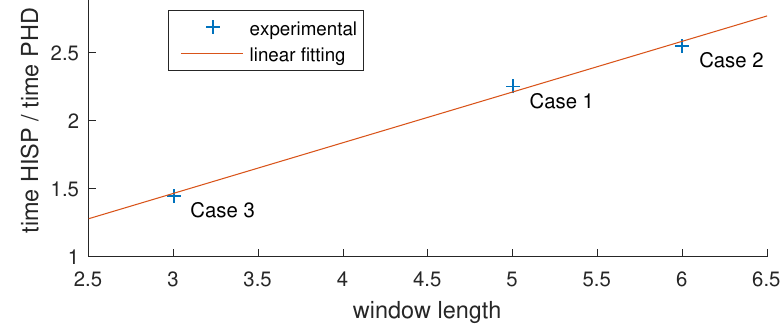}
\caption{Relative computational time of the \gls{hisp} w.r.t.\ the \gls{phd} filter for Cases 1 to 3, averaged over $100$ \gls{mc} runs.}
\label{fig:computationalTime}
\end{figure}

Although both the \gls{phd} and \glspl{hisp} have a linear complexity, the computational cost for the \gls{hisp} tends to be higher than for the \gls{phd} filter, especially when the time window used in track extraction is large. The ratio between the measured computational times of the two filters is displayed in \cref{fig:computationalTime} and shows that the relation between this ratio and the length of the time window appears to be linear, although this is only based on $3$ data points.

\section*{Conclusion}

A new \gls{mtt} algorithm, called the \gls{hisp}, has been derived and detailed. When studying this filter, it appeared that there is more than one way of using the update equations and that there are different possible approximations as well as diverse applicable modelling alternatives. In this sense, the \gls{hisp} can be seen as a general and computationally-affordable way of approaching the problem of \gls{mtt}. The \gls{hisp} allows for characterising each hypothesis separately thus giving a local picture of the underlying multi-target problem while controlling the level of approximation. Its efficiency has been compared with the performance of the \gls{phd}=, \gls{cphd} and \gls{lmb} filters and the results show that the \gls{hisp} outperforms these algorithms in several cases with varying probabilities of detection and statistics for the false alarms.

\bibliographystyle{abbrv}
\bibliography{Thesis}

\appendix

\begin{proof}[Proof of \cref{thm:timeFiltering}]
The principle of the proof is as follows: in general, all combinations of transition functions must be applied to the joint law of all hypotheses at time $t-1$. Without the considered assumptions on the transitions, this would generate a large number of multi-target configurations at time $t$. However, even when making use of the considered modelling, each single-target law which has its probability mass on $\boX^{\inx}$ at time $t-1$ can be predicted using two transitions: $q^{\pi}_t$ and $q^{\omega}_t$. In order to obtain a single predicted multi-target configuration, we introduce a new transition $q_t$ defined as $q_t(\bfx' \given \bfx) = q^{\pi}_t(\bfx' \given \bfx) + q^{\omega}_t(\bfx' \given \bfx)$ for any $\bfx,\bfx' \in \boX$ as well as its extension $\bq_t$ which is defined as in \cref{eq:extensionKernel}. As a consequence, we only consider the transitions $\bq^{\alpha}_t$ and $\bq_t$.

The transition $\bq_t$ is the only one which can be used to predict the law $\delta_{\psi_{\op}}$ without yielding $\delta_{\varphi}$, so that there is no viable alternative. The transition $\bq^{\alpha}_t$ is applied to a fixed number $n^{\alpha}_t$ of targets with law $\delta_{\psi_{\obj}}$ and the other targets with this law can only be propagated with $\bq_t$, so that there is only one possibility. The multiplicity $n_{\varphi}$ of the law $\delta_{\varphi}$ is irrelevant and is not estimated. For the same reasons, targets that are almost-surely in the state space $\boX^{\inx}$ can only be propagated with $\bq_t$ so that the joint law $P^{\obj}_{t|t-1}$ induced by the hypotheses in $\bbI^{\obj}_{t|t-1}$ is
\eqnsml{
P^{\obj}_{t|t-1}(\bfx_1,\dots,\bfx_N) = \prod_{j=1}^{n^{\u}_{t-1}} \bigg[ \int \bq_t(\bfx_j\given\bfx') \bp^{\u}_{t-1}(\bfx')\d\bfx' \bigg] \\
\times \prod_{\bsk \in \bbI^{\m}_{t-1}} \bigg[\int \bq_t(\bfx_{n(\bsi)}\given\bfx') \bp^{\bsk}_{t-1}(\bfx')\d\bfx'\bigg],
}
where $N = n^{\u}_{t-1} + |\bbI^{\m}_{t-1}| = |\bbI^{\obj}_{t|t-1}|$ and where $n$ is any one-to-one mapping between the sets $\bbI^{\m}_{t-1}$ and $\{n^{\u}_{t-1}+1,\dots,N\}$. However, when marginalising all single-target laws but one, say $\bsi \in \bbI^{\m}_{t-1}$, and when considering the restriction to $\boX^{\inx}$ or to $\psi_{\obj}$, the prediction of the single-target law with index $\bsk$ via $\bq^{\pi}_t$ or $\bq^{\omega}_t$ can be recovered through
\eqns{
\bp^{\bsk,\omega}_{t|t-1}(\psi_{\obj}) \defeq \int \delta_{\psi_{\obj}}(\bfx_{n(\bsk)}) P^{\obj}_{t|t-1}(\bfx_1,\dots,\bfx_N) \d(\bfx_1,\dots,\bfx_N)
}
and $\bp^{\bsk,\omega}_{t|t-1}(\varphi) = 1-\bp^{\bsk,\omega}_{t|t-1}(\psi_{\obj})$, and through
\eqns{
\bp^{\bsk,\pi}_{t|t-1} (\bfx) \defeq \int \delta_{\bfx}(\bfx_{n(\bsk)}) P^{\obj}_{t|t-1}(\bfx_1,\dots,\bfx_N) \d(\bfx_1,\dots,\bfx_N)
}
for any $\bfx \in \boX^{\inx}$ and $\bp^{\bsk,\pi}_{t|t-1} = \int_{\boX^{\inx}} \bp^{\bsk,\pi}_{t|t-1} (\bfx) \d\bfx$. It can indeed be verified that
\eqnsa{
\bp^{\bsk,\omega}_{t|t-1} & = \int \bq^{\omega}_t(\cdot \given \bfx) \bp^{\bsk}_{t-1}(\bfx)\d\bfx \\
\bp^{\bsk,\pi}_{t|t-1} & = \int \bq^{\pi}_t(\bfx \given \bfx) \bp^{\bsk}_{t-1}(\bfx)\d\bfx
}
The propagation of the single-target law $\bp^{\u}_{t-1}$ of the never-detected targets via $\bq^{\pi}_t$ or $\bq^{\omega}_t$ can be similarly recovered; for instance the law defined by
\eqns{
\bp^{\u,\pi}_{t|t-1}(\bfx) \defeq \int P^{\obj}_{t|t-1}(\bfx, \bfx_2,\dots,\bfx_N) \d(\bfx_2,\dots,\bfx_N)
}
and by $\bp^{\u,\pi}_{t|t-1}(\varphi) = 1 - \int\bp^{\u,\pi}_{t|t-1}(\bfx)\d\bfx$ is indeed equal to the common law of the never-detected targets predicted via $\bq^{\pi}_t$. Thus, \ref{predApprox} yields a multi-target configuration of the form
\eqnsml{
\{ (\bp^{\bsi}_{t|t-1}, 1)\}_{\bsi \in \bbI^{\m}_{t|t-1}} \cup \{(\bp^{\u,\pi}_{t|t-1}, n^{\u}_{t-1})\} \\ \cup \{(\bp^{\u,\omega}_{t|t-1}, n^{\u}_{t-1})\} \cup \{(\bp^{\alpha}_t, n^{\alpha}_t)\} \cup \bar\calP^{\psi,\varphi}_{t-1}
}
where $\bp^{\alpha}_t = \bq_t^{\alpha}(\cdot \given \psi_{\obj})$. The targets that disappear during prediction are then mixed by Simplification~\ref{it:forgetDisappearedUndetected} and hence return to $\bar\calP^{\psi,\varphi}_{t-1}$. Although the law $\int \bq^{\omega}_t(\cdot \given \bfx) \bp^{\bsk}_{t-1}(\bfx)\d\bfx$ might have some probability mass on $\varphi$, the mixing with an infinite number of laws $\delta_{\psi_{\obj}}$ will make this negligible. Then, by Simplification~\ref{it:mixUndetected}, the law $\bp^{\u}_{t|t-1}$ results from the mixing of the $n^{\u}_{t-1}$ never-detected targets with law $\bp^{\u,\pi}_{t|t-1}$ and the $n^{\alpha}_t$ appearing targets with law $\bp^{\alpha}_t$, so that $n^{\u}_{t|t-1} = n^{\u}_{t-1} + n^{\alpha}_t$ also follows. This concludes the proof of the theorem.
\end{proof}

The proof of \cref{thm:HispObsFilter} requires the introduction of some additional notations. Let $\calA_t$ be the set of subsets of $\bbI^{\obj}_t$ describing all the possible associations of hypotheses in $\bbI^{\obj}_{t|t-1}$ with observations in $\bar{Z}_t$ such that hypotheses do not share non-empty observations and let $\calA_t(\bsk,z)$ be the subset of $\calA_t$ made of sets such that the hypothesis $\bsk \in \bbI_{t|t-1}$ is associated with $z \in \bar{Z}_t$. 
Let $w_t$ be the \gls{disp}'s posterior probability mass function on $\calA_t \times \bbN$ (see \cite[Corollary 3.4]{Houssineau2015}) corresponding to the \gls{hisp}'s predicted multi-target configuration $\bar\calP_{t|t-1}$ and defined for any $I \in \calA_t$ as
\eqnsml{
w_t(I,n^{\u}_{t|t-1} - n_I) = \big[ \bp^{\u,\d}_{t|t-1}(\phi) \big]^{n^{\u}_{t|t-1} - n_I} \bigg[ \prod_{z \in Z^{\op}_I} \bp^{\op,z}_{t|t-1}(z) \bigg] \\
\times \bigg[ \prod_{z \in Z'_t \setminus Z^{\op}_I} \bp^{\op,z}_{t|t-1}(\phi) \bigg] \Bigg[ \prod_{\substack{(\bsk,s,z) \in \sigma_t^{-1}[I] \\ (\bsk,z) \neq (\u,\phi) }} \bp^{\bsk,\d}_{t|t-1}(z) \Bigg],
}%
with $n_I \leq |Z_t|$ the number of targets in $I$ detected for the first time and with $Z^{\op}_I \subseteq Z_t$ the subset of observations considered as false alarms in $I$. The proof of the expression of $w_t$ is out of the scope of the present article; however, the scalar $w^I_t = w_t(I,n^{\u}_{t|t-1} - n_I)$ can be seen to be the marginal likelihood associated with the update of the hypotheses selected by $I$ with the corresponding observations and with the failure of the detection of the others. It follows that
\eqnl{eq:defWt}{
W_t = \sum_{I \in \calA_t} w^I_t
}
is the full multi-target marginal likelihood.

\begin{proof}[Proof of \cref{thm:HispObsFilter}]
For any $\bsi \in \bbI_t$ with $(\bsk,s,z) \defeq \sigma_t^{-1}(\bsi)$, the posterior marginal law $\bp^{\bsi}_t$ on $\bar\boX$ can be characterised by
\eqnl{eq:proofHispBayes1}{
\bp^{\bsi}_t(\bfx) = \sum_{I \in \calA_t(\bsk,z)} \dfrac{w^I_t}{W_t \, \bp^{\bsk,s}_{t|t-1}(z)} L^{s,z}_t(\bfx) p^{\bsk}_{t|t-1}(\bfx),
}
for any $\bfx \in \boX$, and by $\bp^{\bsi}_t(\varphi) = 1 - \int_{\boX} \bp^{\bsi}_t(\bfx) \d\bfx$. 
Defining
\eqns{
w_{\ex}^{\bsk,z} = \sum_{I \in \calA_t(\bsk,z)} w^I_t,
}
we can rewrite \cref{eq:proofHispBayes1} as
\eqns{
\bp^{\bsi}_t(\bfx) = W^{-1}_t w_{\ex}^{\bsk,z} L^{\d,z}_t(\bfx) p^{\bsk}_{t|t-1}(\bfx)
}
for any $\bfx \in \boX$, so that it remains to prove that the marginal likelihood $W_t$ can be equivalently expressed as the denominator of either \cref{eq:perIndividual} or \cref{eq:perObs}. It is sufficient to verify that for any $z \in \bar{Z}_t$ and any $\bsk \in \bbI_{t|t-1}$, it holds that
\eqns{
\bigcup_{\bsk' \in \bbI_{t|t-1}} \calA_t(\bsk',z) = \bigcup_{z' \in \bar{Z}_t} \calA_t(\bsk,z') = \calA_t
}
since the constraint ``hypothesis $\bsk$ is associated with observation $z$'' is loosen by taking the union over all $\bsk$ in $\bbI_{t|t-1}$ or over all $z$ in $Z_t$ (informally, this constraint becomes ``an object is associated with $z$'' or ``$\bsk$ is associated with some possibly-empty observations'', which is always true under the standard \gls{mtt} assumptions). It follows that
\eqns{
W_t = \sum_{I \in \calA_t} w^I_t  = \sum_{\bsk' \in \bbI_{t|t-1}} w_{\ex}^{\bsk',z} \bp^{\bsk'}_t(z)
= \sum_{z' \in \bar{Z}_t} w_{\ex}^{\bsk , z'} \bp^{\bsk}_t(z').
}
The expression \cref{eq:thm:HispObsFilter:fullDistribution} of the approximated multi-target configuration $\bar\calP_t$ if then formed under \ref{upApprox}.
\end{proof}

\Cref{cor:HispFactWeights} is a corollary of the following lemma, where $\calA^{\m}_t$ refers to the subset of $\calA_t$ where all hypotheses have been detected at least once up to time $t-1$.

\begin{lemma}
\label{lem:factoP}
Considering \ref{approx1} for the sets $\bbI^{\m}_t$ and $Z_t$, the multi-target marginal likelihood $W_t$ can be factorised as
\eqns{
W_t = C_t \prod_{\bsk \in \bbI^{\m}_{t|t-1}} \bigg[ \bp^{\bsk,\d}_{t|t-1}(\phi) + \sum_{z \in Z_t} \dfrac{p^{\bsk,\d}_{t|t-1}(z)}{C^{\u,\op}_t(z)} \bigg],
}
where
\eqns{
C_t = \big[\bp^{\u,\d}_{t|t-1}(\phi)\big]^{n^{\u}_{t|t-1}} \bigg[ \prod_{z \in Z'_t} \bp^{\op,z}_{t|t-1}(\phi) \bigg] \bigg[ \prod_{z \in Z_t} C^{\u,\op}_t(z) \bigg].
}
Alternatively, considering \ref{approx2} for $\bbI^{\m}_t$ and $Z_t$, it holds that
\eqns{
W_t = C^{\phi}_t \prod_{z \in Z_t} \bigg[ C^{\u,\op}_t(z) + \sum_{\bsk \in \bbI^{\m}_{t|t-1}} \dfrac{p^{\bsk,\d}_{t|t-1}(z)}{\bp^{\bsk,\d}_{t|t-1}(\phi)} \bigg],
}
where the constant $C^{\phi}_t$ is the joint probability for all the targets in the system to be undetected at time $t$, defined as
\eqns{
C^{\phi}_t = \big[\bp^{\u,\d}_{t|t-1}(\phi)\big]^{n^{\u}_{t|t-1}} \bigg[ \prod_{z \in Z'_t} \bp^{\op,z}_{t|t-1}(\phi) \bigg] \bigg[ \prod_{\bsk \in \bbI^{\m}_{t|t-1}} \bp^{\bsk,\d}_{t|t-1}(\phi) \bigg].
}
\end{lemma}

\begin{proof}
We first rewrite the multi-target marginal likelihood $w^I_t$ in a suitable way, for any $I \in \calA_t$. For this purpose let $J^{\m,\d}_I$ be the subset of $\bbI^{\m}_{t|t-1}$ corresponding to targets with index in $I$ that have been detected at time $t$, let $Z^{\m,\d}_I$ be the subset of $Z_t$ containing the corresponding observations, let $\sigma_I : J^{\m,\d}_I \to Z^{\m,\d}_I$ be the one-to-one mapping describing this identification, and let $Z^{\u}_I$ and $Z^{\op}_I$ be the subsets of $Z_t$ containing the observations associated with the never-detected targets and with false alarms respectively, then it holds that
\eqns{
Z_t = Z^{\m,\d}_I \uplus Z^{\u}_I \uplus Z^{\op}_I,
}
and $w^I_t$ can be expressed as
\eqnsa{
& w^I_t = C^{\phi}_t\\
& 
\!\!\times\!\bigg[ \prod_{(\bsk,z) \in \graph(\sigma_I)} \dfrac{p^{\bsk,\d}_{t|t-1}(z)}{\bp^{\bsk,\d}_{t|t-1}(\phi)} \bigg]\!
\bigg[ \prod_{z \in Z^{\u}_I} \dfrac{p^{\u,\d}_{t|t-1}(z)}{\bp^{\u,\d}_{t|t-1}(\phi)} \bigg]\!
\bigg[ \prod_{z \in Z^{\op}_I} \dfrac{\bp^{\op,z}_{t|t-1}(z)}{\bp^{\op,z}_{t|t-1}(\phi)} \bigg]
}
where $\graph(\sigma_I) \defeq \{(\bsk,\sigma_I(\bsk)) \st \bsk \in J^{\m,\d}_I\}$ is the graph of $\sigma_I$. We can proceed to the second step of the proof by rewriting $W_t$ as follows
\eqnsml{
W_t = C^{\phi}_t \sum_{I \in \calA^{\m}_t} \bigg[ \prod_{(\bsk,z) \in \graph(\sigma_I)} \dfrac{p^{\bsk,\d}_{t|t-1}(z)}{\bp^{\bsk,\d}_{t|t-1}(\phi)} \bigg] \\
\times \Bigg[ \sum_{\substack{Z_{\u}, Z_{\op} \subseteq Z_t \sst\\ Z_{\u} \uplus Z_{\op} = Z_t - Z^{\m,\d}_I}} \bigg[ \prod_{z \in Z_{\u}} \dfrac{p^{\u,\d}_{t|t-1}(z)}{\bp^{\u,\d}_{t|t-1}(\phi)} \bigg] \bigg[ \prod_{z \in Z_{\op}} \dfrac{\bp^{\op,z}_{t|t-1}(z)}{\bp^{\op,z}_{t|t-1}(\phi)} \bigg] \Bigg].
}
We conclude by noticing that the sum over $Z_{\u}$ and $Z_{\op}$ has a binomial form and can thus be factorised, so that
\eqnsa{
W_t & = C^{\phi}_t \sum_{I \in \calA^{\m}_t}\bigg[ \prod_{(\bsk,z) \in \graph(\sigma_I)} \dfrac{p^{\bsk,\d}_{t|t-1}(z)}{\bp^{\bsk,\d}_{t|t-1}(\phi)} \bigg] \bigg[ \prod_{z \in Z_t - Z^{\m,\d}_I} C^{\u,\op}_t(z) \bigg], \\
& = C^{\phi}_t \bigg[ \prod_{z \in Z_t} C^{\u,\op}_t(z) \bigg] \\
& \qquad\qquad \times \sum_{I \in \calA^{\m}_t} \bigg[ \prod_{(\bsk,s,z) \in \sigma_t^{-1}[I]} \dfrac{p^{\bsk,\d}_{t|t-1}(z)}{\bp^{\bsk,\d}_{t|t-1}(\phi) C^{\u,\op}_t(z) } \bigg].
}
We finally notice that
\eqnsa{
& \prod_{\bsk \in \bbI^{\m}_{t|t-1}} \bigg[ \bp^{\bsk,\d}_{t|t-1}(\phi) + \sum_{z \in Z_t} \dfrac{p^{\bsk,\d}_{t|t-1}(z)}{C^{\u,\op}_t(z)} \bigg] \\
& = \! \bigg[ \prod_{\bsk \in \bbI^{\m}_{t|t-1}} \! \bp^{\bsk,\d}_{t|t-1}(\phi) \bigg] \! \sum_{I \in \calA^{\m}_t} \! \bigg[ \prod_{(\bsk,s,z) \in \sigma_t^{-1}[I]} \dfrac{p^{\bsk,\d}_{t|t-1}(z)}{\bp^{\bsk,\d}_{t|t-1}(\phi) C^{\u,\op}_t(z) } \bigg]
}
where the product on the r.h.s.\ has been developed and the obtained terms have been simplified using \ref{approx1}, which completes the proof of the first part of the \lcnamecref{lem:factoP}. The alternative expression of $W_t$ under \ref{approx2} can be proved in a similar way.
\end{proof}


\end{document}

%% file: Kernels3.pdf_tex
\begingroup%
  \makeatletter%
  \providecommand\color[2][]{%
    \errmessage{(Inkscape) Color is used for the text in Inkscape, but the package 'color.sty' is not loaded}%
    \renewcommand\color[2][]{}%
  }%
  \providecommand\transparent[1]{%
    \errmessage{(Inkscape) Transparency is used (non-zero) for the text in Inkscape, but the package 'transparent.sty' is not loaded}%
    \renewcommand\transparent[1]{}%
  }%
  \providecommand\rotatebox[2]{#2}%
  \ifx\svgwidth\undefined%
    \setlength{\unitlength}{270.2824707bp}%
    \ifx\svgscale\undefined%
      \relax%
    \else%
      \setlength{\unitlength}{\unitlength * \real{\svgscale}}%
    \fi%
  \else%
    \setlength{\unitlength}{\svgwidth}%
  \fi%
  \global\let\svgwidth\undefined%
  \global\let\svgscale\undefined%
  \makeatother%
  \begin{picture}(1,0.5144398)%
    \put(0,0){\includegraphics[width=\unitlength,page=1]{Kernels3.pdf}}%
    \put(0.70714232,0.14345533){\color[rgb]{0,0,0}\makebox(0,0)[lb]{\smash{$\psi_{\obj}$}}}%
    \put(0.16289897,0.14845937){\color[rgb]{0,0,0}\makebox(0,0)[lb]{\smash{$\psi_{\op}$}}}%
    \put(0.41,0.36168537){\color[rgb]{0,0,0}\makebox(0,0)[lb]{\smash{$\boX^{\inx}$ at $t-1$}}}%
    \put(-0.00368538,0.16750063){\color[rgb]{0,0,0}\makebox(0,0)[lb]{\smash{$\pi$}}}%
    \put(0.89164394,0.13120686){\color[rgb]{0,0,0}\makebox(0,0)[lb]{\smash{$\pi$}}}%
    \put(0.20828678,0.48070773){\color[rgb]{0,0,0}\makebox(0,0)[lb]{\smash{$\pi$}}}%
    \put(0,0){\includegraphics[width=\unitlength,page=2]{Kernels3.pdf}}%
    \put(0.55843342,0.00923513){\color[rgb]{0,0,0}\makebox(0,0)[lb]{\smash{$\alpha$}}}%
    \put(0,0){\includegraphics[width=\unitlength,page=3]{Kernels3.pdf}}%
    \put(0.73433031,0.2889682){\color[rgb]{0,0,0}\makebox(0,0)[lb]{\smash{$\omega$}}}%
    \put(0,0){\includegraphics[width=\unitlength,page=4]{Kernels3.pdf}}%
    \put(0.32,0.22739832){\color[rgb]{0,0,0}\makebox(0,0)[lb]{\smash{$\boX^{\inx}$ at $t$}}}%
    \put(0,0){\includegraphics[width=\unitlength,page=5]{Kernels3.pdf}}%
  \end{picture}%
\endgroup%